\documentclass[a4j,11pt]{article}
\usepackage{framed}
\usepackage{graphicx}
\newtheorem{theorem}{Theorem}
\newtheorem{proposition}{Proposition}
\newtheorem{corollary}{Corollary}
\newtheorem{lemma}{Lemma}

\def\vec#1{\mbox{\boldmath $#1$}}
\usepackage{amsmath,amssymb}

\usepackage{amssymb}
\usepackage{latexsym}
\def\vec#1{\mbox{\boldmath $#1$}}
\usepackage{amsmath}

\DeclareMathOperator{\diag}{diag}

\begin{document}

\title{Asymptotics for penalized additive $B$-spline regression}

\author{
{\sc Takuma Yoshida}$^{1}$\ \ {\sc and}\ \ {\sc Kanta Naito}$^{2}$\thanks{E-mail:\ {\tt naito@riko.shimane-u.ac.jp}}\\
$^{1}${\it Graduate School of Science and Engineering}\ \ {\it and}\ \ $^{2}${\it Department of Mathematics} \\
{\it Shimane University, Matsue 690-8504, Japan}
}

\date{27th April 2011}

\maketitle

\begin{abstract}
This paper is concerned with asymptotic theory for penalized spline estimator in bivariate additive model. 
The focus of this paper is put upon the penalized spline estimator obtained by the backfitting algorithm. 
The convergence of the algorithm as well as the uniqueness of its solution are shown. 
The asymptotic bias and variance of penalized spline estimator are derived by an efficient use of the asymptotic results for the penalized spline estimator in marginal univariate model. 
Asymptotic normality of estimator is also developed, by which an approximate confidence interval can be obtained.
Some numerical experiments confirming theoretical results are provided. 
\end{abstract}

{\it key words:}{
Additive model;\ Backfitting algorithm;\ $B$-spline;\ Penalized spline.}
\vspace{10mm}

\section{Introduction}
The additive model is a typical regression model with multidimensional covariates and is usually expressed as
$$
y_i=f_1(x_{i1})+\cdots+f_D(x_{iD})+\varepsilon_i,
$$
for given data $\{(y_i,x_{i1},\cdots,x_{iD}):i=1,\cdots,n\}$, where each $f_d (d=1,\cdots,D)$ is a univariate function with a certain degree of smoothness. 
This paper focuses on the bivariate additive model, in which $D=2$.
  
The additive model has become a popular smoothing technique and its fundamental properties have been summarized in literature such as Buja et al. (1989) and Hastie and Tibshirani (1990). 
Buja et al. (1989) proposed the so-called backfitting algorithm, which is efficient for nonparametric estimation of $f_{d}(d=1,\cdots,D)$. 
The backfitting algorithm is a repetition update algorithm and its convergence and the uniqueness of its solution are not always assured.
 Buja et al. (1989) showed the sufficient condition for convergence of the backfitting algorithm and the uniqueness of its solution for the bivariate additive model.

In this paper, we discuss the asymptotic properties of the penalized spline estimator for the additive model with $D=2$. 
Unlike spline smoothing, the asymptotic results of kernel smoothing for the additive model have been obtained. 
Ruppert and Opsomer (1997) showed that a certain kernel smoothing for the additive model satisfies the sufficient condition for convergence of the backfitting algorithm and the uniqueness of its solution. 
Furthermore, they derived the asymptotic bias and variance of the kernel estimator for the bivariate additive model.

Opsomer (2000) presented the sufficient condition for convergence of the backfitting algorithm and the uniqueness of its solution for the $D$-variate additive model in Lemma 2.1. 
The asymptotic bias and variance of the kernel estimator for the $D$-variate additive model were also derived under the assumption that the sufficient condition for convergence of the backfitting algorithm holds. 
Wand (1999) investigated  asymptotic normality of the kernel estimator for the $D$-variate additive model by elegant use of the results in Opsomer (2000).
We observe from Wand's results of asymptotic normality that kernel estimators of $f_d$'s are asymptotically independent.

Many researchers have explored the effectiveness of spline smoothing, such as Wahba (1975) and Green and Silverman (1994).
Penalized spline estimators have been discussed in O'Sullivan (1986), Eilers and Marx (1996), Marx and Eilers (1998) and Ruppert et al. (2003). 
Despite its richness of application, asymptotics for spline smoothing seems have not yet been sufficiently developed. 

For the univariate model ($D=1$), Agarwal and Studen (1980) and Zhou et al. (1998) obtained important asymptotic results for the regression spline. 
Hall and Opsomer (2005) gave the mean squared error and consistency of the penalized spline estimator. 
The asymptotic bias and variance of the penalized spline estimator were obtained in Claeskens et al. (2009). 
Kauermann et al. (2009) worked with the generalized linear model.  
Wang et al. (2011) showed that the penalized spline estimator is asymptotically equivalent to a Nadaraya-Watson estimator. 
Thus, it seems that developments of asymptotic theories of the penalized spline are relatively recent events and we note that those works are mainly regarding the univariate model ($D=1$). 
In the case of multidimensional covariates, Stone (1985) showed the consistency of the regression spline in the $D$-variate additive model, but it is not penalized spline.

The aim of this paper is to derive asymptotic bias, asymptotic variance, and asymptotic distribution of the penalized spline estimator in the bivariate additive model.  
The penalized spline estimator for the bivariate additive model is obtained using the penalized least squares method and the backfitting algorithm. 
The uniqueness of the solution of the backfitting algorithm cannot be proved in general, but its convergence property can be shown.
However, it is demonstrated that the solution of the backfitting algorithm is asymptotically unique and the objective function for the penalized least squares method is shown to be asymptotically convex.  
As will be seen in the subsequent section, the penalized spline estimator in a bivariate setting has a closed form, which we can use for asymptotic manipulations.
The properties of band matrices play an important role as a mathematical tool in asymptotic considerations.  
The effect of the initial value required for implementing the backfitting algorithm is also investigated.

This paper is organized as follows. 
In Section 2, our model setting and estimating equation in the penalized least squares method are discussed and the backfitting algorithm to obtain the solution is composed. 
Section 3 provides the asymptotic bias and variance of the penalized spline estimator and then its asymptotic normality is developed. 
Furthermore, the uniqueness of the solution of the backfitting algorithm is discussed. 
Section 4 includes numerical studies to validate the theory and an application to real data is reported. 
In Section 5, some suggestions that are necessary to develop the asymptotics for the general $D$-variate spline additive model are noted by comparing similar results already developed for the kernel estimator. 
Proofs for theoretical results are all given in the Appendix.

\section{Model setting}

\subsection{Bivariate additive spline model}

Consider a bivariate additive regression model
\begin{eqnarray}
y_i=f_1(x_{i1})+f_2(x_{i2})+\varepsilon_i, \label{model1}
\end{eqnarray}
for data $\{(y_i,x_{i1},x_{i2}):i=1,\cdots,n\}$, 
where $f_j(\cdot)$ is an unknown regression function and $\varepsilon_i$'s are independent random errors with $E[\varepsilon_i]=0$ and $V[\varepsilon_i]=\sigma^2(x_{i1},x_{i2})<\infty$. 
We assume $E[f_j(X_j)]=0(j=1,2)$ to ensure identifiability of $f_j$.
Let $q_j(x_j)$ be the density of $X_j$ and $q(x_1,x_2)$ be the joint density of $(X_1,X_2)$. 
We assume without loss of generality that $(x_{i1},x_{i2})\in(0,1)\times (0,1)$ for all $i\in \{1,\cdots,n\}$.   
  
Now we consider the $B$-spline model
$$
s_j(x_{j})=\sum_{k=-p+1}^{K_n} B_k^{[p]}(x_{j})b_{j,k}
$$ 
as an approximation to $f_j(x_j)$ at any $x_{j} \in (0,1)$ for $j=1,2$. 
Here, 
$B_k^{[p]}(x)(k=-p+1,\cdots,K_n)$ are $p$th degree $B$-spline basis functions defined recursively as 
\begin{eqnarray*}
B_k^{[0]}(x)&=&
\left\{
\begin{array}{cc}
1,& \kappa_{k-1}<x\leq \kappa_k,\\
0,& {\rm otherwise},
\end{array}
\right. \\
B_k^{[p]}(x)&=&\frac{x-\kappa_{k-1}}{\kappa_{k+p-1}-\kappa_{k-1}}B_k^{[p-1]}(x)+\frac{\kappa_{k+p}-x}{\kappa_{k+p}-\kappa_{k}}B_{k+1}^{[p-1]}(x) ,
\end{eqnarray*}
where $\kappa_k=k/K_n (k=-p+1,\cdots,K_n+p)$ are knots, $K_n=O(n^\gamma)$ with $0<\gamma<1/2$, 
and $b_{j,k}( j=1,2, k=-p+1,\cdots,K_n)$ are unknown parameters. 
We denote $B_k^{[p]}(x)$ as $B_k(x)$ in what follows since only the $p$th degree is treated. 
The details and many properties of the $B$-spline function are clarified in de Boor (2001). 
We aim to obtain an estimator of $f_j$ via the $B$-spline additive regression model 
\begin{eqnarray}
y_i=s_1(x_{i1})+s_2(x_{i2})+\varepsilon_i, \label{Bspmodel}
\end{eqnarray}
instead of the model (\ref{model1}).  
The model (\ref{Bspmodel}) can be expressed as 
$$
\vec{y}=X_1\vec{b}_1+X_2\vec{b}_2+\vec{\varepsilon}
$$
by using the notations  
$\vec{y}=(y_1\ \cdots\ y_n)^\prime$, $\vec{b}_1=(b_{1,-p+1} \cdots\ b_{1,K_n})^\prime$, $\vec{b}_2=(b_{2,-p+1} \cdots\ b_{2,K_n})^\prime$, $X_1=(B_{-p+j}(x_{i1}))_{ij}$, $X_2=(B_{-p+j}(x_{i2}))_{ij}$ and $\vec{\varepsilon}=(\varepsilon_1\ \cdots\ \varepsilon_n)^\prime$.
We adopt the estimators $(\hat{\vec{b}}_{1}',\hat{\vec{b}}_{2}')$ of $(\vec{b}_{1}',\vec{b}_{2}')$ defined as the minimizer of
\begin{eqnarray}
L(\vec{b}_1,\vec{b}_2)=(\vec{y}-X_1\vec{b}_1-X_2\vec{b}_2)^\prime(\vec{y}-X_1\vec{b}_1-X_2\vec{b}_2)+\sum_{j=1}^2\lambda_{jn}\vec{b}_j^\prime Q_m\vec{b}_j, \label{pen}
\end{eqnarray}
where $\lambda_{jn} (j=1,2)$ are smoothing parameters and $Q_m$ is the $m$th order difference matrix.  
This criterion is called the penalized least squares method and it has been frequently utilized in spline regression (Eilers and Marx (1996)). 
For a fixed point $x_j\in(0,1)$, the estimator $\hat{f}_j(x_j)$ of $f_j(x_j)$
is
$$
\hat{f}_j(x_j)=\sum_{k=-p+1}^{K_n} B_k(x_{j})\hat{b}_{j,k}.
$$ 
and is called the penalized spline estimator of $f_j(x_j)$. 
The predictor of $y$ at a fixed point $(x_1,x_2)\in(0,1)\times (0,1)$ is defined as
$$
\hat{y}=\hat{f}_1(x_1)+\hat{f}_2(x_2).
$$
 
Since $E[f_j(X_j)]=0$ is assumed for $f_j$, the estimator of each component $f_j$ is usually centered. 
Hence $\hat{f}_j(x_j)$ is rewritten as
$$
\hat{f}_{j,c}(x_j)=\hat{f}_j(x_j)-\frac{1}{n}\sum_{i=1}^n \hat{f}_j(x_{ij}),
$$ 
as discussed in Wang and Yang (2007).
In this paper, however, we do not examine $\hat{f}_{j,c}$ because our interests are in asymptotics for $\hat{f}_j$ and $\hat{y}$, 
and asymptotic distributions of $\hat{f}_j(x_j)$ and $\hat{f}_{j,c}(x_j)$ become equivalent.

\subsection{Backfitting algorithm}\label{bfit}

Let $\vec{b}=(\vec{b}_1^\prime\ \vec{b}_2^\prime)^\prime$. 
In general, $\hat{\vec{b}}=(\hat{\vec{b}}_1^\prime\ \hat{\vec{b}}_2^\prime)^\prime$ is a solution of 
\begin{eqnarray}
\frac{\partial L(\vec{b}_1,\vec{b}_2)}{\partial \vec{b}}=\vec{0}. \label{kai}
\end{eqnarray}
In fact, the solution of (\ref{kai}) can be written as 
$\vec{b}_1=\Lambda_1^{-1}X_1^\prime(\vec{y}-X_2\vec{b}_2)$ and $\vec{b}_2=\Lambda_2^{-1}X_2^\prime(\vec{y}-X_1\vec{b}_1), \label{backz}
$
where $\Lambda_j=X_j^\prime X_j+\lambda_{jn}Q_m$.
However, this method has one defect:\ the $L(\vec{b}_1,\vec{b}_2)$ is not in general convex as the function of $\vec{b}$. 
Hence, the solution of (\ref{kai}) does not necessarily become the minimizer of (\ref{pen}). 
Marx and Eilers (1998) also noted this point as a typical problem of additive spline regression. 

Let $\tilde{\vec{b}}=(\tilde{\vec{b}}_1'\ \tilde{\vec{b}}_2')'$ be a minimizer of (\ref{pen}). 
Then it is important to investigate the difference between $\hat{\vec{b}}$ and $\tilde{\vec{b}}$ asymptotically. 
If the difference is vanishingly small, it shows that $\hat{\vec{b}}$ asymptotically minimizes (\ref{pen}). 
The details of this assertion are given in Section \ref{minim}. 

In this paper, our estimator of $(\vec{b}_{1}', \vec{b}_{2}')'$ is composed by using the backfitting algorithm obtained from the solution of (\ref{kai}). 
The merit and usage of the backfitting algorithm are clarified in Hastie and Tibshirani (1990).   
The $\ell$-stage backfitting estimators $\vec{b}_1^{(\ell)}$ and $\vec{b}_2^{(\ell)}$ are defined as 
$$
\vec{b}_1^{(\ell)}=\Lambda_1^{-1}X_1^\prime(\vec{y}-X_2\vec{b}_2^{(\ell -1)})\ \ \ \mbox{and}\ \ \ \vec{b}_2^{(\ell)}=\Lambda_2^{-1}X_2^\prime(\vec{y}-X_1\vec{b}_1^{(\ell)}),
$$
respectively, where $\vec{b}_2^{(0)}$ is an initial value. 
Then, the $\ell$-stage backfitting estimator $f_j^{(\ell)}(x_j)$ of $f_j(x_j)$ at $x_j\in(0,1)$ is obtained as 
$$
f_j^{(\ell)}(x_j)=\sum_{k=-p+1}^{K_n} B_k(x_{j})b^{(\ell)}_{j,k}=\vec{B}(x_j)^\prime\vec{b}^{(\ell)}_j,\ \ j=1,2,
$$
where $\vec{B}(x_j)=(B_{-p+1}(x_j)\ \cdots\ B_{K_n}(x_j))^\prime$. 
A mathematical property of the backfitting algorithm is that $\vec{b}^{(\infty)}=(\vec{b}_1^{(\infty)'},\vec{b}_2^{(\infty)'})'\equiv \lim_{\ell\rightarrow \infty}(\vec{b}_1^{(\ell)'}, \vec{b}_2^{(\ell)'})'$ satisfies 
\begin{eqnarray}
\frac{\partial L(\vec{b})}{\partial \vec{b}}\Big|_{\vec{b}=\vec{b}^{(\infty)}}=\vec{0}. \label{kai2}
\end{eqnarray}
The backfitting algorithm itself is applicable in not only bivariate but also the general $D$-variate additive model. 
However, $\vec{b}_j^{(\ell)}$ can be explicitly expressed only for the case $D=2$.
By referring to (5.24) on page 119 of Hastie and Tibshirani (1990), $\vec{b}_j^{(\ell)}$ can be calculated as  
\begin{eqnarray}
\label{backf}
\left.
\begin{array}{lll}
\vec{b}^{(\ell)}_1
&=&\displaystyle (X_1^\prime X_1)^{-1}X_1^\prime\vec{y}-(X_1^\prime X_1)^{-1}X_1^\prime\sum_{j=0}^{\ell-1} \{S_1S_2\}^j(I_n-S_1)\vec{y}\\
&&
-(X_1^\prime X_1)^{-1}X_1^\prime(S_1S_2)^{\ell-1}S_1X_2\vec{b}_2^{(0)}, \\
\vec{b}^{(\ell)}_2&=&\displaystyle \Lambda_2^{-1}X_2^\prime\sum_{j=0}^{\ell-1} \{S_1 S_2\}^j(I_n-S_1)\vec{y}
+\Lambda_2^{-1}X_2^\prime(S_1 S_2)^{\ell-1} S_1X_2\vec{b}_2^{(0)},
\end{array}
\right.
\end{eqnarray}
where $S_j=X_j\Lambda_j^{-1}X_j^{\prime}$. 
It is shown by Theorem 10 of Buja et al. (1989) that $\vec{b}_j^{(\infty)}(j=1,2)$ converge depending on $\vec{b}_2^{(0)}$. 
Thus, the backfitting estimators $\vec{b}^{(\infty)}_{1}$ and $\vec{b}^{(\infty)}_{2}$ converge, but the vectors to which they converge are not unique, depending on the initial value. 
We will study the asymptotic behavior of $f_j^{(\infty)}(x_j)=\vec{B}(x_j)^\prime \vec{b}^{(\infty)}_j$, as well as the relationship of $\vec{b}^{(\infty)}$ and $\tilde{\vec{b}}$ from now on.

\section{Asymptotic theory}
We prepare some symbols and notations to be used hereafter. 
Let $I_n$ be the identity matrix of size $n$.
Define a matrix $G_k=(G_{k,ij})_{ij}$, where the $(i,j)$-component is 
$$
G_{k,ij}=\int_0^1 B_i(x)B_j(x)q_k(x)dx
$$
for $k=1,2$.
Define a matrix $\Sigma_{k}=(\Sigma_{k,ij})_{ij}$, where the $(i,j)$-component is
$$
\Sigma_{k,ij}=\int_0^1 \int_0^1 \sigma^2(x_1,x_2)B_i(x_{k})B_j(x_{k})q(x_1,x_2)dx_1dx_2
$$
for $k=1,2$.

Let a vector $\vec{b}_j^*$ be such that $\vec{B}(\cdot)\vec{b}^*_j$ satisfies the best $L_{\infty}$ approximation to the true function $f_j$. 
For further information on this point, see Zhou et al. (1998).

For a matrix $X_n=(X_{ij,n})_{ij}$, if $\displaystyle\max_{i,j}\{n^{\alpha}|X_{ij,n}|\}=O_{P}(1)(o_{P}(1))$, then it is written as $X_n=O_{P}(n^{-\alpha}\vec{1}\vec{1}^\prime)(o_{P}(n^{-\alpha}\vec{1}\vec{1}^\prime))$. 
This notation will be used for matrices with fixed sizes and sizes depending on $n$.
 
In spline smoothing, the smoothing parameter $\lambda_{jn}$ is usually selected as $\lambda_{jn}\rightarrow \infty$ with $n\rightarrow \infty$ because 
a spline curve often yields overfitting for large $n$. 
In the following, we assume that $\lambda_{jn}=o(nK_n^{-1})$. 
Hence, we choose $\lambda_{jn}$ as $\lambda_{jn}\rightarrow \infty$ and $\lambda_{jn}=o(nK_n^{-1})$.

\subsection{Asymptotic distribution of the penalized spline estimator}\label{assplline}
Let $f^{(\ell)}_{0j}(x_j)$ be $f^{(\ell)}_{j}(x_j)$ with $\vec{b}_2^{(0)}=\vec{0}$. 
Then, $f^{(\ell)}_{1}(x_1)$ and $f^{(\ell)}_{2}(x_2)$ with arbitrary initial value $\vec{b}_2^{(0)}$ can be expressed as
\begin{eqnarray*}
f_1^{(\ell)}(x_1)=f_{01}^{(\ell)}(x_1)-\vec{B}(x_1)^\prime (X_1^\prime X_1)^{-1}X_1^\prime(S_1S_2)^{\ell-1}S_1X_2\vec{b}_2^{(0)}
\end{eqnarray*}
and 
\begin{eqnarray*}
f_2^{(\ell)}(x_2)=f_{02}^{(\ell)}(x_2)+\vec{B}(x_2)^\prime \Lambda_2^{-1}X_2(S_1 S_2)^{\ell-1} S_1X_2\vec{b}_2^{(0)},
\end{eqnarray*}
respectively.
First, we investigate the influence of $\vec{b}_2^{(0)}$ on $f^{(\ell)}_{j}(x_j)$, which is summarized as follows.

\vspace{5mm}

\begin{proposition}\label{initial}
Suppose that $\lambda_{jn}=o(nK_n^{-1})$.
Then for $\ell=1,2,\cdots,$
\begin{eqnarray*}
|f_1^{(\ell)}(x_1)-f_{01}^{(\ell)}(x_1)|=O_P(K_n^{-2\ell+1})\ \ \ \mbox{and}\ \ \ |f_2^{(\ell)}(x_2)-f_{02}^{(\ell)}(x_2)|=o_P(K_n^{-2\ell+1}),\ \ as\ \  n\rightarrow \infty.
\end{eqnarray*}
In particular, as $\ell\rightarrow \infty$, 
$
f_j^{(\infty)}(x_j)=f_{0j}^{(\infty)}(x_j),\ \ j=1,2.
$
\end{proposition}

\vspace{5mm}

Proposition \ref{initial} claims that the influence of $\vec{b}_2^{(0)}$ on $f_{j}^{(\ell)}$ can be ignored for large $n$. 
In other words, for any initial value $\vec{b}_2^{(0)}\in\mathbb{R}^{K_n+p}$, we can uniquely obtain $f_j^{(\infty)}(x_j)$ as $n\rightarrow \infty$.
Hence, it suffices to consider $f_{0j}^{(\ell)}(x_j)$ instead of $f_{j}^{(\ell)}(x_j)$ to develop asymptotics under manipulations $\ell\rightarrow \infty$ and $n\rightarrow \infty$.
Here, $f_{01}^{(\ell)}(x_1)$ and $f_{02}^{(\ell)}(x_2)$ can be written as 
$$
f_{01}^{(\ell)}(x_1)=f_{01}^{(1)}(x_1)-\vec{B}(x_1)^{\prime}(X_1^\prime X_1)^{-1}X_1^\prime\sum_{k=1}^{\ell-1} \{S_1S_2\}^k(I_n-S_1)\vec{y}
$$
and 
$$
f_{02}^{(\ell)}(x_2)= f_{02}^{(1)}(x_2)+ \vec{B}(x_2)^{\prime}\Lambda_2^{-1}X_2^\prime\sum_{k=1}^{\ell-1} \{S_1 S_2\}^k(I_n-S_1)\vec{y},
$$
respectively. 
This allows the following.

\vspace{5mm}

\begin{proposition}\label{backfit}
Suppose that $\lambda_{jn}=o(nK_n^{-1})$. 
Then for any fixed point $(x_1,x_2)\in(0,1)\times (0,1)$,
\begin{eqnarray*}
f_{01}^{(\infty)}(x_1)=f_{01}^{(1)}(x_1)+O_P(K_n^{-1})\ \ \ \mbox{and}\ \ \ f_{02}^{(\infty)}(x_2)=f_{02}^{(1)}(x_2)+o_P(K_n^{-1}),\ as\ n\rightarrow \infty.
\end{eqnarray*}
\end{proposition}

\vspace{5mm}

We see from (\ref{backf}) that 
\begin{eqnarray*}
f_{01}^{(1)}(x_1)=\vec{B}(x_1)^\prime \Lambda_1^{-1}X_1^\prime \vec{y}\ \ \ \mbox{and}\ \ \ f_{02}^{(1)}(x_2)=\vec{B}(x_2)^\prime \Lambda_2^{-1}X_2^\prime (I_n-S_1)\vec{y}.
\end{eqnarray*}
It should be noted that $f_{01}^{(1)}(x_1)$ has the same form as the penalized spline estimator based on the dataset $\{(y_i,x_{i1}):i=1,\cdots,n\}$ in the univariate regression model ($D=1$). 
This form is very important because the asymptotic bias and variance of the penalized spline estimator for the univariate regression model have been already derived by Claeskens et al. (2009). 
Similarly, $f_{02}^{(1)}(x_2)$ includes $\vec{B}(x_2)^\prime \Lambda_2^{-1}X_2^\prime\vec{y}$, which is the same as the penalized spline estimator for univariate regression based on $\{(y_i,x_{i2}):i=1,\cdots,n\}$.

We denote $f_{j}^{(\infty)}(x_j)$ as $\hat{f}_j(x_j)=\vec{B}(x_j)^\prime\hat{\vec{b}}_j$, which does not depend on the initial value $\vec{b}_2^{(0)}$ as $n\rightarrow \infty$.
The usefulness of Propositions \ref{initial} and \ref{backfit} is that we can realize the asymptotic equivalence between the backfitting estimator $\hat{f}_{j}$ and the (marginal) univariate penalized spline estimator.  
By using the results of Claeskens et al. (2009), we obtain Theorem 1.

\vspace{5mm}

\begin{theorem}\label{mv}
Suppose that $f_j\in C^{p+1}$ and $\lambda_{jn}=o(nK_n^{-1})$. 
Then for any fixed point $(x_1,x_2)\in(0,1)\times (0,1)$, as $n\rightarrow \infty$, 
\begin{eqnarray*}
E[\hat{f}_j(x_j)]&=&f_j(x_j)+b_{j,\lambda}(x_j)+O_P(K_n^{-1})+o_P(\lambda_{jn}K_nn^{-1}),\\
V\left[
\begin{array}{c}
\hat{f}_1(x_1)\\
\hat{f}_2(x_2)
\end{array}
\right]
&=&
\left[
\begin{array}{cc}
V_{1}(x_1)(1+o_p(1))&O_P(n^{-1})\\
O_P(n^{-1})&V_{2}(x_2)(1+o_P(1))
\end{array}
\right], 
\end{eqnarray*}
where, 
\begin{eqnarray*}
b_{j,\lambda}(x_j)=-\frac{\lambda_{jn}}{n}\vec{B}(x)^\prime G_j^{-1}Q_m\vec{b}_{j}^*,\ \ V_j(x_j)=\frac{1}{n}\vec{B}(x_j)^\prime G_{j}^{-1}\Sigma_{j}G_{j}^{-1}\vec{B}(x_j).
\end{eqnarray*}
\end{theorem}

\vspace{5mm}

By using Theorem \ref{mv}, we have the asymptotic joint distribution of $[\hat{f}_1(x_1)\ \hat{f}_2(x_2)]^\prime$.

\vspace{5mm}

\begin{theorem}\label{clt}
Suppose that there exists $\delta>0$ such that $E[|\varepsilon_i|^{2+\delta}]<\infty$ and $f_j\in C^{p+1}$. 
Furthermore, $\gamma$ and $\lambda_{jn}$ satisfy $1/3<\gamma<1/2$ and $\lambda_{jn}=o((nK_n^{-1})^{1/2})$. 
Then for any fixed point $(x_1,x_2)\in(0,1)\times (0,1)$, as $n\rightarrow \infty$, 
\begin{eqnarray*}
V\left[
\begin{array}{c}
\hat{f}_1(x_1)\\
\hat{f}_2(x_2)
\end{array}
\right]
^{-1/2}
\left[
\begin{array}{c}
\hat{f}_1(x_1)-f_1(x_1)\\
\hat{f}_2(x_2)-f_2(x_2)
\end{array}
\right]
\xrightarrow {D}N_2\left( 
\left[
\begin{array}{c}
0\\ 
0 
\end{array} 
\right]
,I_2\right).
\end{eqnarray*}
\end{theorem}

\vspace{5mm}

From Theorem \ref{clt}, $\hat{f}_1(x_1)$ and $\hat{f}_2(x_2)$ are asymptotically independent. 
Asymptotic normality and the independence of $\hat{f}_1(x_1)$ and $\hat{f}_2(x_2)$ in kernel smoothing also hold, as shown in Wand (1999). 
Thus, the penalized spline estimator and the kernel estimator for the additive model have the same asymptotic property. 
Asymptotic normality of $\hat{y}$ can be shown as a direct consequence of Theorem \ref{clt}. 
We briefly note the pointwise confidence interval for $f_j(x_j)$ by exploiting the distribution of $\hat{f}_j(x_j)$ obtained from Theorem \ref{clt}.
Here, we treat $\sigma^2(x_{i1},x_{i2})$ as known for all $i\in\{1,\cdots,n\}$, but it should be estimated in data analysis.

\vspace{5mm}

\begin{corollary}\label{confidence}
A $100(1-\alpha)\%$ asymptotic confidence interval of $f_j(x_j)$ at any fixed point $x_j\in(0,1)$ is
$$
\hat{f}_j(x_j)\pm z_{\alpha/2}\sqrt{V[\hat{f}_j(x_j)]},
$$
where $z_{\alpha/2}$ is the $(1-\alpha/2)$th normal percentile. 
\end{corollary}

\vspace{5mm}

The confidence interval in Corollary \ref{confidence} will be applied to a set of real data in Section 4, in which we need to prepare an estimate of $V[\hat{f}_j(x_j)]$.

\subsection{Minimizer of $L(\vec{b})$}\label{minim}
Here, we discuss the difference between $\hat{\vec{b}}=\vec{b}^{(\infty)}$ and $\tilde{\vec{b}}$, the minimizer of (\ref{pen}). 
Although $\hat{\vec{b}}$ is the solution of (\ref{kai}), the problem of whether it minimizes $L(\vec{b})$ or not is not trivial. 
That is, many solutions of (\ref{kai}) might exist because $L(\vec{b})$ is not convex. 
Let the $M$ solutions of (\ref{kai}) be $\{\tilde{\vec{b}}^1,\cdots,\tilde{\vec{b}}^M\}$. 
Then, for any $\vec{b}_2^{(0)}\in\mathbb{R}^{K_n+p}$, there exists $m\in\{1,\cdots,M\}$ such that $\hat{\vec{b}}=\vec{b}^{(\infty)}=\tilde{\vec{b}}^{m}$.   
However, $\vec{b}^{(\infty)}$ is asymptotically not dependent on $\vec{b}_2^{(0)}$ as implied in Proposition \ref{initial}. 
Therefore, the uniqueness of the penalized spline estimator obtained by the backfitting algorithm is asymptotically satisfied. 
Furthermore, Theorem 3 says that $\hat{\vec{b}}$ minimizes $L(\vec{b})$.

\vspace{5mm}

\begin{theorem}\label{convex}
Let $H(L)$ be the Hessian matrix of $L(\vec{b})$. 
Then  
$H(L)$ is asymptotically positive definite.
\end{theorem}

\vspace{5mm}

We see that asymptotic properties of the penalized spline estimator for the additive model can be obtained not only by Theorems \ref{mv} and \ref{clt}, but also by Theorem \ref{convex}.

\section{Numerical studies}

In this section, we see the behavior of the estimator and validate Theorem 2 numerically by simulation. 
In addition, we aim to obtain an asymptotic confidence interval using a real dataset. 
We utilize the cubic spline ($p=3$) and the second order difference penalty ($m=2$) in all of the following numerical studies.

\subsection{Simulation}
We choose the true functions $f_1(x_1)=\sin(2\pi x_1)$, $f_2(x_2)=2^{-1}\cos(\pi x_2)$ and the error is $\varepsilon_i\sim U(-0.5,0.5)$. 
Here, $U(a,b)$ is a uniform distribution on an interval $[a,b]$. 
The explanatory variables $x_{ij} (i=1,\cdots,n,j=1,2)$ are derived from $x_{ij}\sim U(0,1)$. 
Then, $f_1$ and $f_2$ satisfy $E[f_1(X_1)]=0$ and $E[f_2(X_2)]=0$, respectively.
We demonstrate three simulations. 

In Simulation-1, we compare $f^{(\ell)}_{j}(x_j)$ with the true $f_j(x_j)$. 

In Simulation-2, we compare $f_{j}^{(\ell)}(x_j)$ with
$$
\hat{f}_{pen,j}(x_j)=\vec{B}(x_j)^\prime \Lambda_j^{-1}X_j^\prime \vec{y},
$$
which is the penalized spline estimator for univariate regression based on $(y_i,x_{ij})$.

In Simulation-3, we compare the density of $N_2(\vec{0},I_2)$ with the kernel density estimate of simulated
\begin{eqnarray}
V\left[
\begin{array}{c}
f^{(\ell)}_1(x_1)\\
f^{(\ell)}_2(x_2)
\end{array}
\right]
^{-1/2}
\left[
\begin{array}{c}
f^{(\ell)}_1(x_1)-f_1(x_1)\\
f^{(\ell)}_2(x_2)-f_2(x_2)
\end{array}
\right] \label{2norm}
\end{eqnarray}
to validate Theorem \ref{clt}, where we note that the covariance matrix of $[f^{(\ell)}_1(x_1)\ f^{(\ell)}_2(x_2)]'$ can be exactly calculated and it in fact was used in this simulation. 
The bandwidth of the kernel density estimate is selected by the method of Sheather and Jones (1991).
The algorithm of Simulation-3 is given as follows:

\vspace{5mm}

\begin{enumerate}
\item[Step 1]\ Generate $x_{ij}\sim U(0,1)$ for $j=1,2,i=1,\cdots,n$.
\item[Step 2]\ Generate the data $\{(y_{i},x_{i1},x_{i2})|i=1,\cdots,n\}$ from (\ref{model1}) and $\varepsilon_i\sim U(-0.5,0.5)$.
\item[Step 3]\ Calculate $f^{(\ell)}_j(x_j)$ at fixed point $x_1=x_2=0.5$.
\item[Step 4]\ Calculate the values of (\ref{2norm}).
\item[Step 5]\ Iterate from Step 2 to Step 4, 1000 times.
\item[Step 6]\ Draw the kernel density estimate of (\ref{2norm}) and compare with the density of $N_2(\vec{0},I_2)$.
\end{enumerate}

\vspace{5mm}

The results of Simulation-1, Simulation-2 and Simulation-3 are displayed in Figure 1, Figure 2, and Figure 3, respectively. 
In all simulation settings, $K_n=2n^{2/5}$, $\lambda_{1n}=\lambda_{2n}=2n^{2/5}K_n^{-1/2}$ and $\ell=10$ were adopted. 
We set the sample size $n=1000$ for Simulation-1 and Simulation-2, and $n=100$ and $n=1000$ for Simulation-3.

\begin{figure}
\begin{center}
\includegraphics[width=60mm,height=55mm]{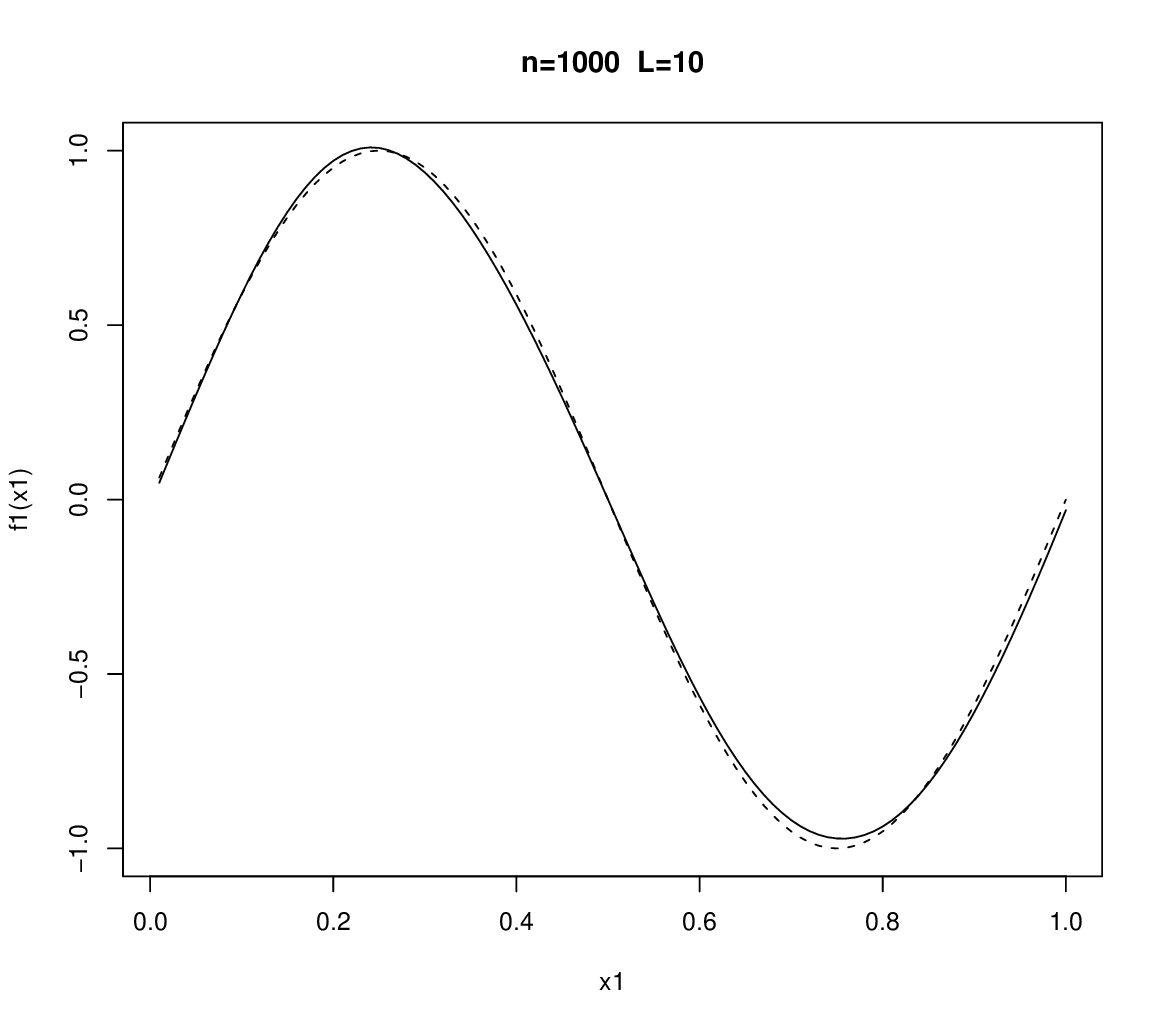}
\includegraphics[width=60mm,height=55mm]{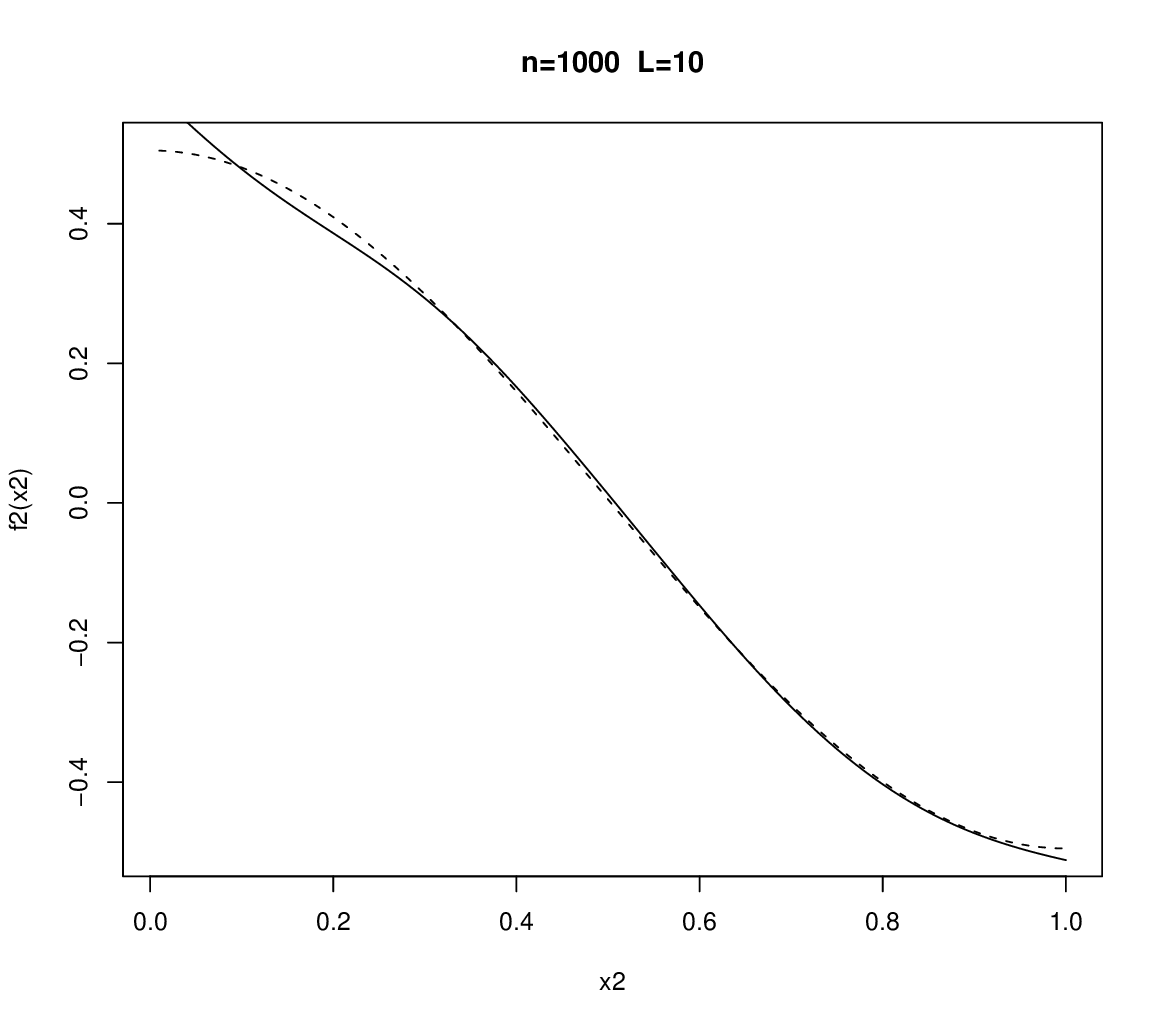}\\
\caption{The curves of $f_j^{(\ell)}(x_j)$ and $f_j(x_j)$.
The left panel is for $f_1^{(\ell)}$(solid line) and $f_1$(dashed line).  
The right panel is for $f_2^{(\ell)}$(solid) and $f_2$(dashed). }
\end{center}
\end{figure}

We see from Figure 1 that the backfitting estimator $f_j^{(10)}$ approximates $f_j$ well. 
We also observe in Figure 2 that the differences between $f^{(10)}_{j}$ and $\hat{f}_{pen,j}(x_j)$ are small, which means that $\hat{f}_{pen,j} \approx f^{(1)}_{j}$ dominates the backfitting estimator as claimed in Proposition \ref{backfit}.

The contour plots of the density estimate of (\ref{2norm}) and of the density of $N_2(\vec{0},I_2)$ are drawn in Figure 3. 
We observe that there is still a gap between the density estimate and the density of $N_2(\vec{0},I_2)$ in $n=100$. 
However, we see from the case $n=1000$ that the density estimate is clearly approaching the density of $N_2(\vec{0},I_2)$, as claimed in Theorem \ref{clt}.

\begin{figure}
\begin{center}
\includegraphics[width=60mm,height=60mm]{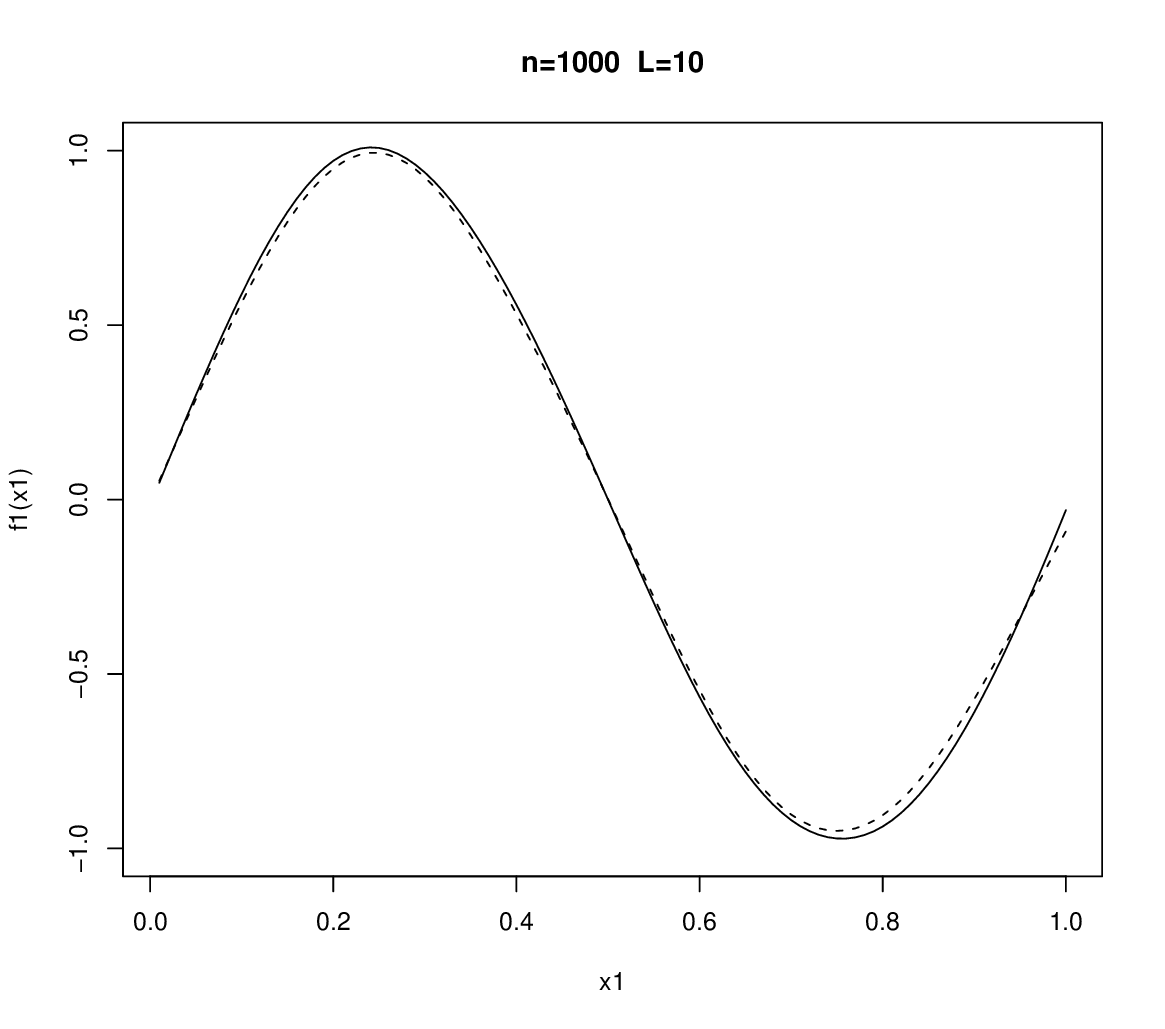}
\includegraphics[width=60mm,height=60mm]{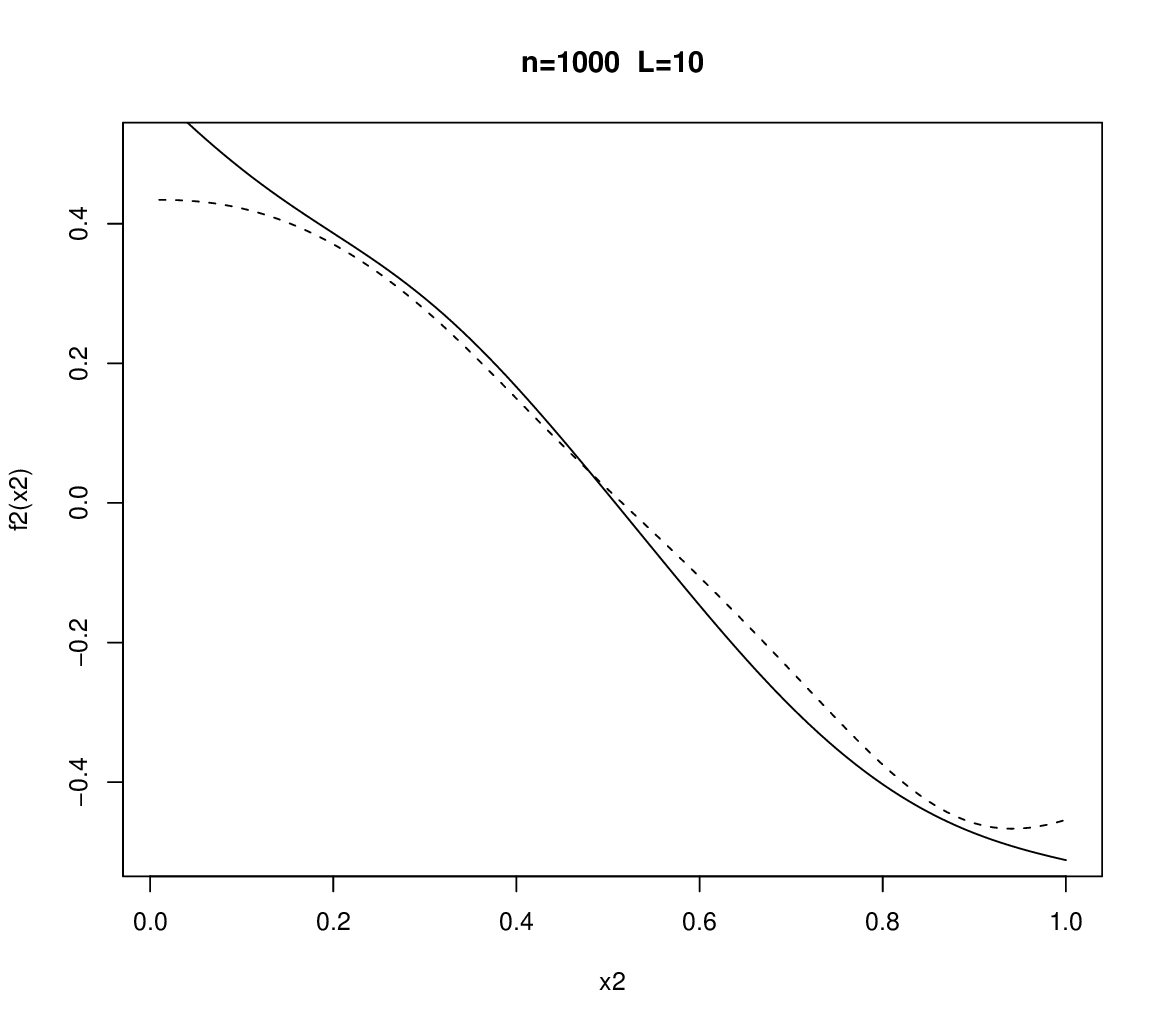}\\
\caption{
The curve of $f_j^{(\ell)}(x_j)$ and $\hat{f}_{pen,j}(x_j)$. The lines $f_1^{(\ell)}(x_1)$ (solid) and $\hat{f}_{pen,1}(x_1)$ (dashed) are drawn in the left panel. 
The lines $f_2^{(\ell)}(x_2)$ (solid) and the $\hat{f}_{pen,2}(x_2)$ (dashed) are drawn in the right panel. }
\end{center}
\end{figure}

\begin{figure}
\begin{center}
\includegraphics[width=60mm,height=60mm]{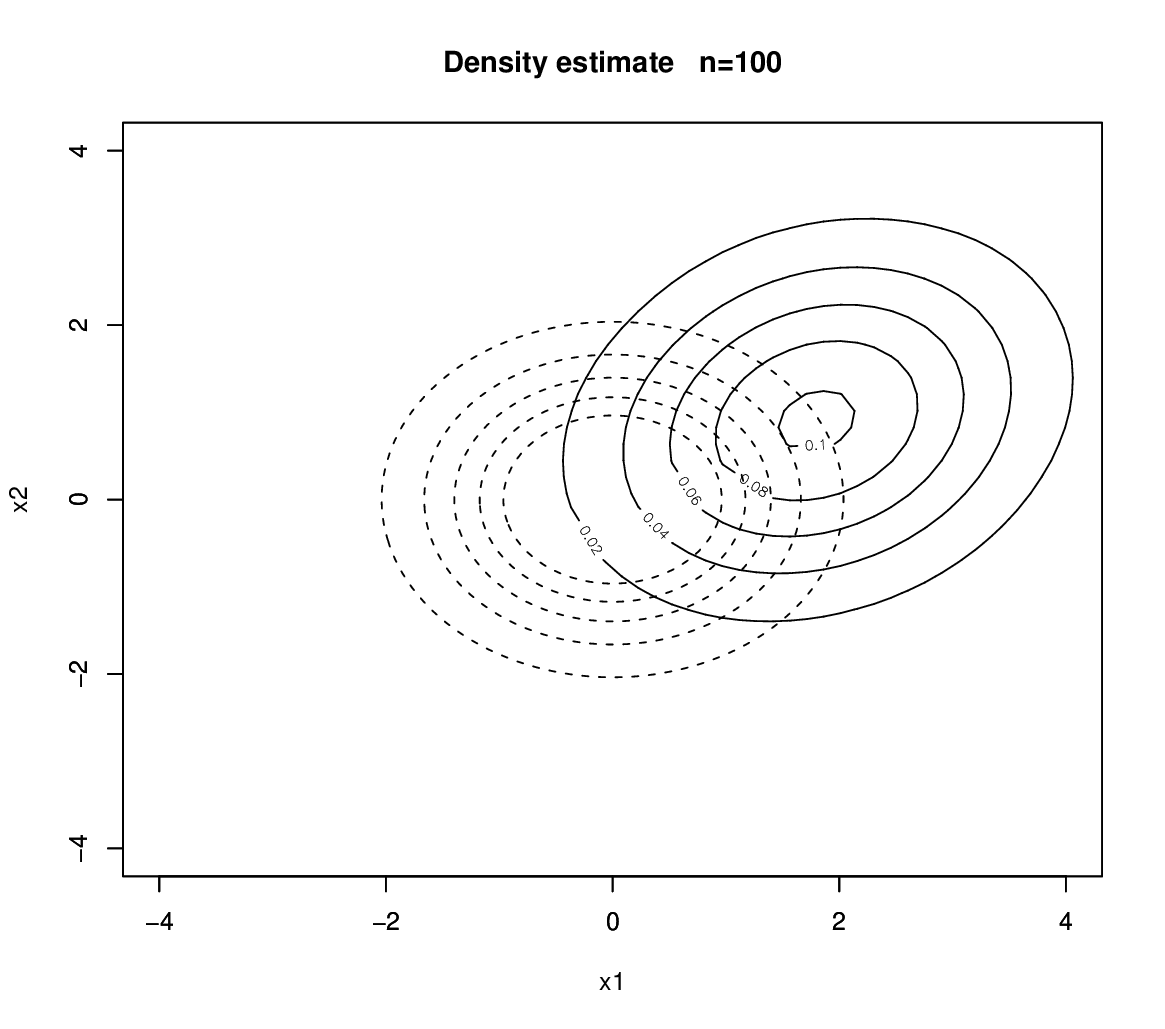}
\includegraphics[width=60mm,height=60mm]{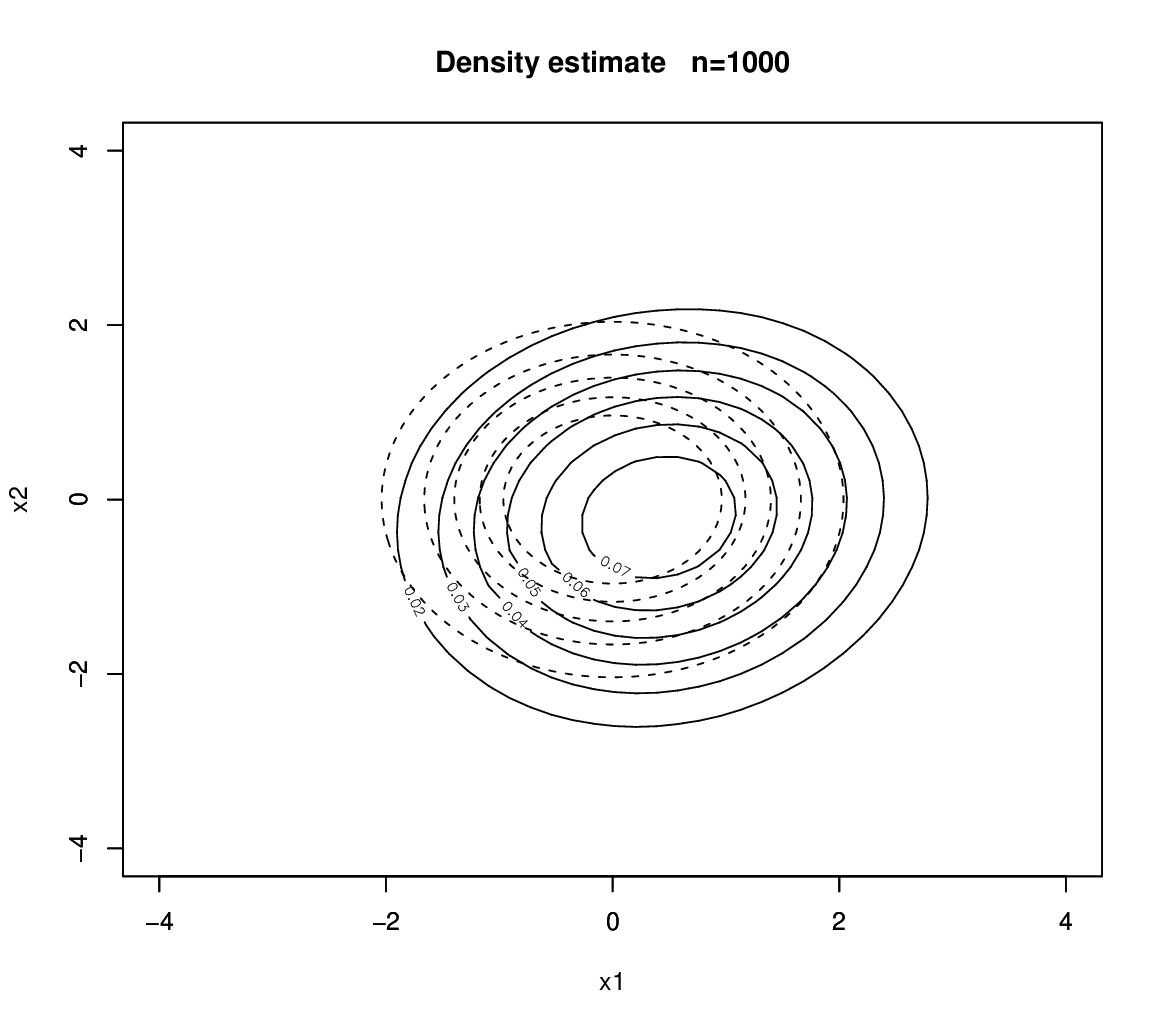}\\
\caption{
The density estimate of (\ref{2norm}) (solid line) and the density of $N(\vec{0},I_2)$ (dashed line). 
The left panel is $n=100$, and the right panel is  $n=1000$. 
The contour lines of $N(\vec{0},I_2)$ are 0.02, 0.04, 0.06, 0.08 and 0.1.}
\end{center}
\end{figure}

\subsection{Application to real data}

We construct the asymptotic pointwise confidence interval of $f_j(x_j)$ by using real data. 
We utilize ozone data with $n=111$ (Hastie et al. (2001)).
We use model (\ref{model1}), where $y$ is ozone concentration (ppb), $x_1$ is daily maximum temperature (${}^\circ\mathrm{C}$) and $x_2$ is wind speed (mph). 
Each $y_i$ is centered and $x_{ij}$'s are modified as $x_{ij}/\max_{1\leq i\leq n}{x_{ij}}$. 
We composed the backfitting estimator $f^{(\ell)}_j(x_j)$ and asymptotic pointwise confidence interval of $f_j(x_j)$ under the assumption that $\sigma^2(x_1,x_2)=\sigma^2$, which can be estimated by
$$
\hat{\sigma}^2=\frac{1}{n}\sum_{i=1}^n \{y_i-f_1^{(\ell)}(x_{i1})-f^{(\ell)}_2(x_{i2})\}^2.
$$  
Again, we used $K_n=2n^{2/5}$, $\lambda_{1n}=\lambda_{2n}=2n^{2/5}K_n^{-1/2}$ and $\ell=10$.

Hastie and Tibshirani (1990) estimated $f_j(j=1,2)$ by using a pseudo additive method based on a smoothing spline. 
In addition, they constructed a pointwise error bar defined as $\hat{f}_j\pm 2\times{\rm standard\ error}$, which is drawn in Figure 9.9 of Hastie and Tibshirani (1990). 
The asymptotic pointwise confidence interval exhibited in Figure 4 looks quite similar to the error bar. 
However, we see that, the asymptotic intervals given in Figure 4 are both smoother than the error bars.
Although this is only an application to one dataset, we thus confirm that the confidence intervals based on asymptotic normality can be applied to real data.

\begin{figure}
\begin{center}
\includegraphics[width=60mm,height=55mm]{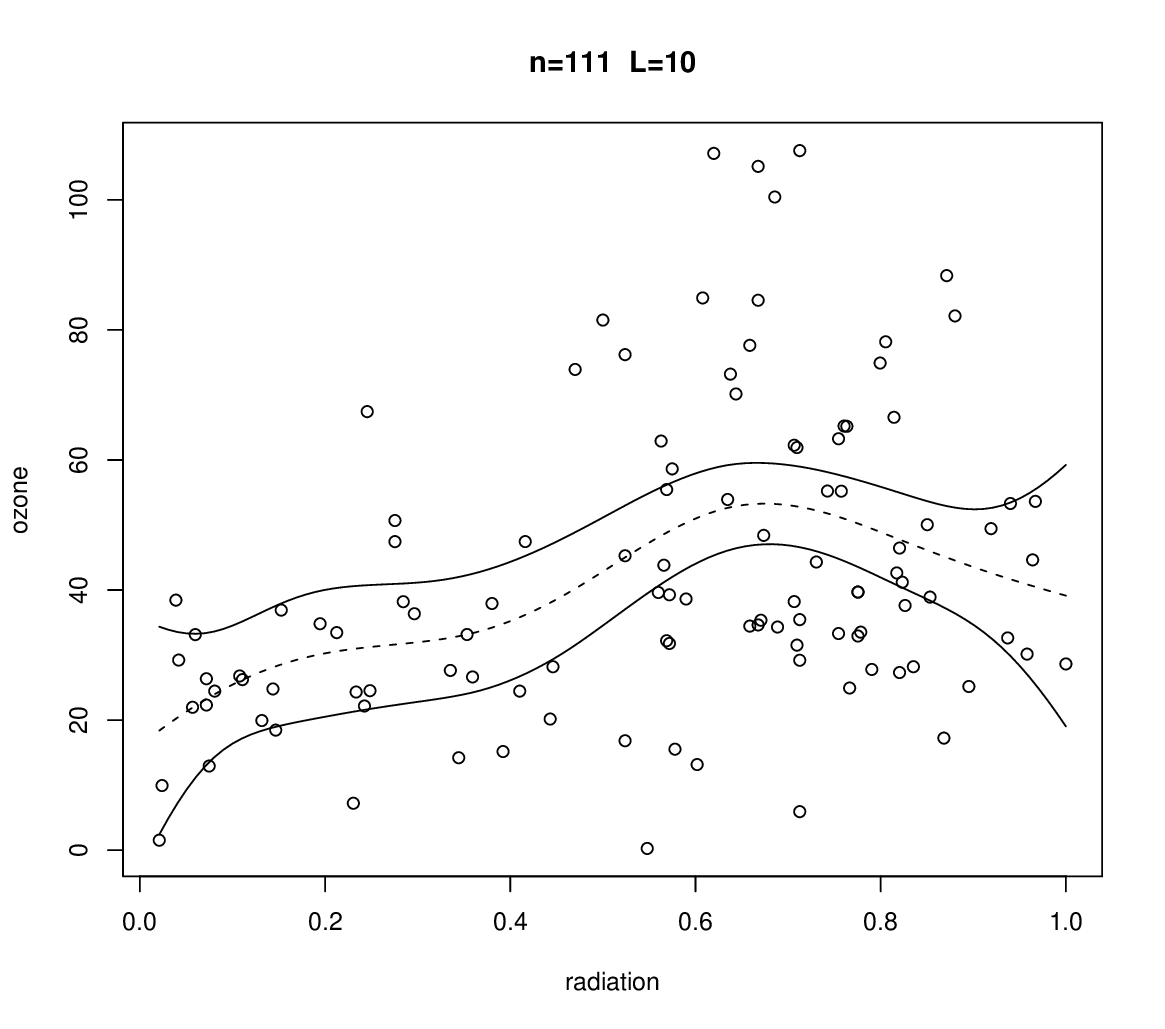}
\includegraphics[width=60mm,height=55mm]{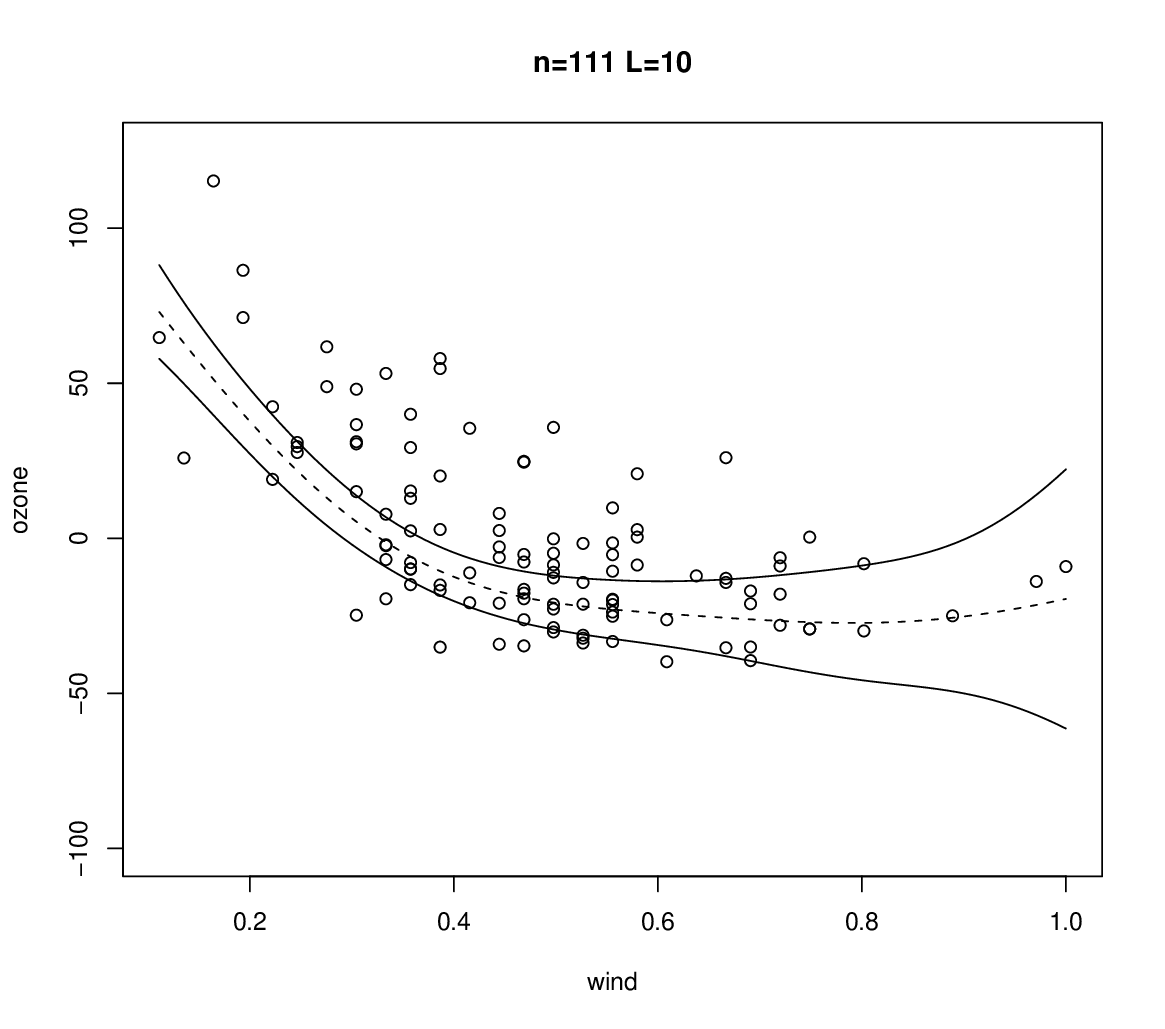}\\
\caption{The asymptotic pointwise confidence intervals for $f_j(x_j)$: the left panel is $(y_i,x_{i1})$ and the right is $(y_i,x_{i2})$.
 The solid lines are $95\%$ confidence intervals and the dashed line is $f_j^{(\ell)}(x_j)$ in both panels.}
\end{center}
\end{figure}

\section{Discussion}
In this paper, asymptotic behavior of the penalized spline estimator in the bivariate additive model is investigated. 
The research in this paper can be seen as a spline version of the work by Ruppert and Opsomer (1997) and Wand (2000). 
To consider a generalization of the work in this paper to the $D$-variate additive model, it might be worthwhile to review the work by Opsomer (2000), including local polynomial fitting in the $D$-variate additive model, as introduced in Section 1.
Let $\vec{f}_d=(f_d(x_{1d})\ \cdots\ f_d(x_{nd}))^\prime (d=1,\cdots,D)$.
Then a formal estimating equation yields the estimator $\hat{\vec{f}}_d$ of $\vec{f}_d$ as
\begin{eqnarray}
\left[
\begin{array}{c}
\hat{\vec{f}}_1\\
\hat{\vec{f}}_2\\
\vdots\\
\hat{\vec{f}}_D
\end{array}
\right]
=
\left[
\begin{array}{cccc}
I_n&S_1&\cdots&S_1\\
S_2&I_n&\cdots&S_2\\
\vdots&\vdots&\ddots&\vdots\\
S_D&S_D&\cdots&I_n
\end{array}
\right]
^{-1}
\left[
\begin{array}{c}
S_1\\
S_2\\
\vdots\\
S_D
\end{array}
\right]
\vec{y}
\equiv 
M^{-1}C\vec{y}\label{eq2}
\end{eqnarray}
provided that $M^{-1}$ exists, where $S_{d}(d=1,...,D)$ are kernel smoothers, as discussed in Opsomer (2000). 
In practice, the estimator is composed by the backfitting algorithm 
\begin{eqnarray}
\label{backkernel}
\left.
\begin{array}{lll}
\vec{f}^{(\ell)}_1&=&S_1(\vec{y}-\vec{f}^{(\ell-1)}_2 \cdots -\vec{f}^{(\ell-1)}_{D}),\\
&\vdots&\\
\vec{f}^{(\ell)}_D&=&S_D(\vec{y}-\vec{f}^{(\ell)}_2 \cdots -\vec{f}^{(\ell)}_{D-1})\\
\end{array}
\right.
\end{eqnarray}
instead of (\ref{eq2}) for reformation of the computational efficiency of $M^{-1}$. 
If $M^{-1}$ exists, it is known that the $\vec{f}^{(\ell)}_d$ in (\ref{backkernel}) converges to the unique $\hat{\vec{f}}_d$ in (\ref{eq2}) as $\ell\rightarrow \infty$.
Opsomer (2000) assumes the sufficient condition for the existence of $M^{-1}$, by which the asymptotic bias and variance of the backfitting estimator for the $D$-variate additive model can be obtained.
It is shown by Ruppert and Opsomer (1997) that $M^{-1}$ certainly exists for the case $D=2$. 
Thus, we see that even in kernel smoothing, such a generalization from bivariate to $D$-variate in the additive model includes the mathematical difficulty.

On the other hand, in the spline method for $D>2$, the smoother is $S_d=X_d\Lambda_d^{-1}X_d^\prime$ and $\vec{f}_d=S_d\vec{b}_d$, where $\vec{b}_d$ is an unknown parameter vector. 
The corresponding matrix $M$ does not have the inverse, even for $D=2$, as detailed in Marx and Eilers (1998). 
Thus, the estimator of $\vec{f}_d$ cannot be written in the form of (\ref{eq2}) and so it might not be reasonable to assume the existence of $M^{-1}$ as the kernel method did. 
The reason why we could proceed with asymptotics for $f^{(\infty)}_j$ is that the explicit form of the backfitting estimator $f^{(\infty)}_j$ can be obtained,
which seems to be impossible for the case $D\geq 3$.
Currently, the only result in this paper that can be generalized to $D\geq 3$ is Theorem 3.

Although it is beyond the scope of this paper, it might be possible to discuss the asymptotics for the penalized spline in the generalized additive model (GAM) in a similar manner. 
Kauermann et al. (2009) studied asymptotic properties of spline regression in the univariate generalized linear model. 
Therefore, the asymptotic theory of the penalized spline in the GAM may be considered for further research.

\section*{Appendix}
For the proofs of Propositions 1-2 and Theorems 1-3, we define $G_{jn}=n^{-1}X_j^\prime X_j(j=1,2)$, $G_{12n}=n^{-1}X_1^\prime X_2$, $G_{21n}=G_{12n}^\prime$, and $\Lambda_{jn}=n^{-1}\Lambda_j(j=1,2)$. 
We need additional lemmas as follows.

\vspace{5mm}

\begin{lemma}\label{G12}
$G_{jn}$, $G_{12n}$ and $\Lambda_{jn}$ satisfy $G_{jn}=O_P(K_n^{-1}\vec{1}\vec{1}^\prime)$, $G_{12n}=O_P(K_n^{-2}\vec{1}\vec{1}^{-1})$ and $\Lambda^{-1}_{jn}=O_P(K_n\vec{1}\vec{1}^\prime)$. 
\end{lemma}
\noindent
{\it Proof of Lemma \ref{G12}}:\ For $j=1,2$, proofs for $G_{jn}=G_j+o_P(K_n^{-1}\vec{1}\vec{1}^\prime)$, $G_j=O(K_n^{-1}\vec{1}\vec{1}^\prime)$ and $\Lambda^{-1}_{jn}=O_P(K_n\vec{1}\vec{1}^\prime)$ have been already given in Claeskens et al. (2009). 
Hence we are going to show $G_{12n}=O_P(K_n^{-2}\vec{1}\vec{1}^{\prime})$. 

Let $G_{12n}=(g_{ij,n})_{ij}$. 
The $(k,h)$-component of $G_{12n}$ is 
$$
g_{kh,n}=\frac{1}{n}\sum_{i=1}^n B_{-p+k}(x_{i1})B_{-p+h}(x_{i2}).
$$
Then $g_{kh,n}$ can be asymptotically expressed as
\begin{eqnarray*}
g_{kh,n}=E[B_{-p+k}(X_1)B_{-p+h}(X_2)]+O_P\left(\sqrt{  \frac{1}{n}V[B_{-p+k}(X_{1})B_{-p+h}(X_{2})] }\right).
\end{eqnarray*}
The $E[B_{-p+k}(X_1)B_{-p+h}(X_2)]$ is bounded by 
\begin{eqnarray*}
&&\min_{u,v\in(0,1)}\{q(u,v)\}\int_0^1B_{-p+k}(x)dx\int_0^1B_{-p+h}(y)dy\\
&&\leq
\int_0^1\int_0^1 B_{-p+k}(x)B_{-p+h}(y)q(x,y)dxdy\\
&&\leq \max_{u,v\in(0,1)}\{q(u,v)\}\int_0^1B_{-p+k}(x)dx\int_0^1B_{-p+h}(y)dy.
\end{eqnarray*}
Hence we get 
$E[B_{-p+k}(X_1)B_{-p+h}(X_2)]=O(K_n^{-2})$ because $\int_0^1B_{-p+k}(x)dx=K_n^{-1}$, see de Boor (2001). 
Similarly, 
$$
E[\{B_{-p+k}(X_{1})B_{-p+h}(X_{2})\}^2]=\int_0^1\int_0^1 \{B_{-p+k}(x)\}^2\{B_{-p+h}(y)\}^2q(x,y)dxdy=O(K_n^{-2}).
$$
Hence we have 
\begin{eqnarray*}
\frac{1}{n}V[B_{-p+k}(X_{1})B_{-p+h}(X_{2})]=O(K_n^{-2}n^{-1}).
\end{eqnarray*}
Since $n^{-1/2}=o(K_n^{-1})$, we have $g_{kh,n}=O_P(K_n^{-2})$. $\Box$

\vspace{5mm}

\noindent
The $G_{jn}(j=1,2)$ and $G_{12n}$ are band matrices:\ for the $(i,k)$-component of $G_{jn}$ and $G_{12n}$, if $|i-k|\leq p$, it is positive and it is 0 if $|i-k|>p$.

\vspace{5mm}

\begin{lemma}\label{band}
Let $A=(a_{ij})_{ij}$ and $B=(b_{ij})_{ij}$ be $K_n\times K_n$ matrices. 
Assume that $K_n\rightarrow \infty$ as $n\rightarrow \infty$, $A=O_P(n^\alpha\vec{1}\vec{1}^\prime)$ and $B$ has $b_{ij}=0$ if $|i-j|>p$ and $b_{ij}=O_P(n^\beta)$ if $|i-j|\leq p$, where $\alpha,\beta\in\mathbb{R}$.
Then 
$AB=O_P(n^{\alpha+\beta}\vec{1}\vec{1}^\prime)$.
\end{lemma}
\noindent
{\it Proof of Lemma \ref{band}}:\ By structural assumption of $B$, the $(i,j)$-component of $AB$ is 
$$
\sum_{k=1}^{K_n} a_{ik} b_{kj}=\sum_{j-p\leq k\leq j+p} a_{ik} b_{kj}=O_P(n^{\alpha+\beta}). \ \ \ \Box
$$

\vspace{5mm}

\begin{lemma}\label{Demko}
Let $A=(a_{ij})_{ij}$ and $B=(b_{ij})_{ij}$ be $K_n\times K_n$ matrices. 
Assume that $K_n\rightarrow \infty$ as $n\rightarrow \infty$, $A=O_P(n^\alpha\vec{1}\vec{1}^\prime)$ and there exist constants $C>0$, $\mu\in(0,1)$ such that $|b_{ij}|\leq n^\beta C\mu^{|i-j|}(1+o_P(1))$ for $1\leq i,j\leq K_n$.
Then 
$AB=O_P(n^{\alpha+\beta}\vec{1}\vec{1}^\prime)$.
\end{lemma}
\noindent
{\it Proof of Lemma \ref{Demko}}:\ The $(i,j)$-component of $AB$ can be evaluated as 
\begin{eqnarray*}
|(AB)_{ij}|
=
\left|\sum_{k=1}^{K_n} a_{ik} b_{kj}\right|\leq 
\max_{1\leq k\leq K_n}\{|a_{ik}|\} \sum_{k=1}^{K_n}|b_{kj}|.
\end{eqnarray*}
Since we have
\begin{eqnarray*}
\sum_{k=1}^{K_n}|b_{kj}|
&\leq &n^\beta C(1+o_P(1)) \sum_{k=1}^{K_n} \mu^{|k-j|}\\
&=& n^\beta C(1+o_P(1))\left(\sum_{k=1}^{j}\mu^{j-k} + \sum_{k=j+1}^{K_n} \mu^{k}\right)\\
&\leq &n^\beta C(1+o_P(1))\left(\sum_{k=0}^{\infty}\mu^{k} + \sum_{k=0}^{\infty} \mu^{k}\right)\\
&=&n^\beta C(1+o_P(1))\frac{2}{1-\mu},
\end{eqnarray*}
it follows that 
$$
|(AB)_{ij}|\leq \max_{1\leq k\leq K_n}\{|a_{ik}|\}n^\beta C\frac{2}{1-\mu}(1+o_P(1))
$$
and hence $AB=O_P(n^{\alpha+\beta}\vec{1}\vec{1}^\prime)$. $\Box$

\vspace{5mm}

\begin{lemma}\label{cla2}
Suppose that $\lambda_{jn}=o(nK_n^{-1})$.
Then for $j=1,2$, there exist constants $C_j>0$ and $\mu_j\in(0,1)$ such that $|(\Lambda^{-1}_{jn})_{ik}|\leq K_nC_j\mu_j^{|i-k|}(1+o_P(1))$, where $(\Lambda^{-1}_{jn})_{ik}$ is the $(i,k)$-component of $\Lambda_{jn}^{-1}$.
\end{lemma}
\noindent
{\it Proof of Lemma \ref{cla2}}:\ Let $\Lambda_{j\cdot}=G_j+\lambda_{jn}n^{-1}Q_m$.
The $\Lambda_{jn}$ can be written as 
\begin{eqnarray*}
\Lambda_{jn}
&=&
G_{jn}+\frac{\lambda_{jn}}{n}Q_m\\
&=&
G_{j}+\frac{\lambda_{jn}}{n}Q_m+G_{jn}-G_j\\
&=&
\Lambda_{j\cdot}+(G_{jn}-G_j)
\end{eqnarray*}
by Lemma \ref{G12}. 
Hence we have 
\begin{eqnarray*}
\Lambda_{jn}^{-1}
&=&
\left(\Lambda_{j\cdot}+(G_{jn}-G_j)\right)^{-1}\\
&=&
\Lambda_{j\cdot}^{-1}\left(I_{K_n+p}+(G_{jn}-G_j)\Lambda_{j\cdot}\right)^{-1}.
\end{eqnarray*}
Let $A_n=-(G_{jn}-G_j)\Lambda_{j\cdot}^{-1}$. 
By Lemma \ref{eigen}, the maximum eigenvalue of $G_{jn}-G_j$ becomes $o_P(K_n^{-1})$ and hence $||(G_{jn}-G_j)||_{2}\leq o_P(K_n^{-1})$. 
Here, for $n\times n$ matrix $A$, $||A||_{2}=\displaystyle\sup_{\vec{x}\not=\vec{0}}\{||A \vec{x}||/||\vec{x}||\}$, where $||\vec{x}||=\sqrt{\vec{x}^\prime\vec{x}}$ for $\vec{x}\in\mathbb{R}^n$.
Further $||\Lambda_{j\cdot}^{-1}||_{2}\leq O(K_n)$ also can be obtained by the proof of Lemma \ref{Demko} and Lemma A1 in Claeskens et al. (2009).
Hence there exists $N_{0}$ such that for any $n>N_{0}$, $||A_n||_{2}\leq ||(G_{jn}-G_j)||_{2} ||\Lambda_{j\cdot}^{-1}||_{2}<1$.
From well-known result of matrix theory, for $n>N_{0}$, $(I_{K_n+p}-A_n)^{-1}$ exists and equality
\begin{eqnarray*}
(I_{K_n+p}-A_n)^{-1}=\sum_{k=0}^{\infty}A_n^k=I_{K_n+p}+A_n \left(\sum_{k=0}^{\infty}A_n^k\right)
\end{eqnarray*}
holds.
So we have  
$$
\sum_{k=0}^{\infty}A_n^k=O_P(\vec{1}\vec{1}^\prime).
$$
Lemmas \ref{G12}, \ref{band} and \ref{Demko} yield
\begin{eqnarray*}
A_n \left(\sum_{k=0}^{\infty}A_n^k\right)
&=&
-(G_{jn}-G_j)\Lambda_{j\cdot}^{-1}\left(\sum_{k=0}^{\infty}A_n^k\right)\\
&=&
-(G_{jn}-G_j)\Lambda_{j\cdot}^{-1}O_P(\vec{1}\vec{1}^\prime)\\
&=&
-(G_{jn}-G_j)O_P(K_n\vec{1}\vec{1}^\prime)\\
&=&
o_P(\vec{1}\vec{1}^\prime)
\end{eqnarray*}
for $n>N_{0}$. 
Therefore the $\Lambda_{jn}^{-1}$ can be asymptotically expressed as 
$$
\Lambda_{jn}^{-1}=\Lambda_{j\cdot}^{-1}\left(I_{K_n+p}+(G_{jn}-G_j)\Lambda_{j\cdot}^{-1}\right)^{-1}=\Lambda_{j\cdot}^{-1}\{ I_{K_n+p}+o_P(\vec{1}\vec{1}^\prime) \}
$$
and its $(i,k)$-component of $\Lambda_{jn}^{-1}$ becomes 
$$
(\Lambda_{jn}^{-1})_{ik}=(\Lambda_{j\cdot}^{-1})_{ik}+ o_P(K_n)
$$
because $\Lambda_{j\cdot}^{-1}o_P(\vec{1}\vec{1}^\prime)=o_P(K_n\vec{1}\vec{1}^\prime)$ by Lemma \ref{Demko}. 
From Claeskens et al. (2009), there exist constants $C_j>0$ and $\mu_j\in(0,1)$ such that $|(\Lambda_{j\cdot}^{-1})_{ik}|\leq K_nC_j\mu^{|i-k|}$.
Hence we finally have  
\begin{eqnarray*}
|(\Lambda_{jn}^{-1})_{ik}|
&\leq&
|(\Lambda_{j\cdot}^{-1})_{ik}|+o_P(K_n)\\
&\leq&
K_nC_j\mu_j^{|i-k|}+o_P(K_n)\\
&=&K_nC_j\mu_j^{|i-k|}\{ 1+o_P(1) \}. \Box
\end{eqnarray*}

\vspace{5mm}

\begin{lemma}\label{S1S2}
For $\ell\geq 1$, there exists a matrix $R^{(\ell)}_{n}=O_P(K_n^{-2(\ell-1)}\vec{1}\vec{1}^\prime)$ such that $\{S_1S_2\}^{\ell}=n^{-1}X_1R^{(\ell)}_{n}X_2^\prime$. 
\end{lemma}
\noindent
{\it Proof of Lemma \ref{S1S2}}:\ We use the inductive method.
First we have the expression 
$$
S_{1}S_{2}=\frac{1}{n^2}X_1\Lambda^{-1}_{1n}X_1^\prime X_2\Lambda_{2n}^{-1}X_2^\prime=\frac{1}{n}X_1\Lambda^{-1}_{1n}G_{12n}\Lambda^{-1}_{2n}X_2^\prime. 
$$
Let $R^{(1)}_{n}=\Lambda^{-1}_{1n}G_{12n}\Lambda^{-1}_{2n}$. 
Then by Lemmas \ref{G12}, \ref{band}, \ref{Demko} and \ref{cla2}, we have $R^{(1)}_{n}=O_P(\vec{1}\vec{1}^\prime)$. 
Next we assume that 
the $(S_1S_2)^\ell$ can be expressed as  
$$
(S_1S_2)^\ell=\frac{1}{n}X_1 R^{(\ell)}_{n}X_2^\prime,
$$
where $R^{(\ell)}_{n}=O_P(K_n^{-2(\ell -1)}\vec{1}\vec{1}^\prime)$. 
For $\ell+1$, 
\begin{eqnarray*}
(S_1S_2)^{\ell+1}
&=&
\frac{1}{n}X_1 R^{(\ell)}_{n}X_2^\prime S_1S_2\\
&=&
\frac{1}{n^2}X_1 R^{(\ell)}_{n}X_2^\prime X_1\Lambda^{-1}_{1n}G_{12n}\Lambda^{-1}_{2n}X_2^\prime\\
&=&
\frac{1}{n}X_1 R^{(\ell)}_{n}G_{21n}\Lambda^{-1}_{1n}G_{12n}\Lambda^{-1}_{2n}X_2^\prime.
\end{eqnarray*}
So we shall put $R^{(\ell+1)}_{n}=R^{(\ell)}_{n}G_{21n}\Lambda^{-1}_{1n}G_{12n}\Lambda^{-1}_{2n}$. 
From Lemma \ref{band}, we get $R^{(\ell)}_{n}G_{21n}=O_P(K_n^{-2\ell}\vec{1}\vec{1}^\prime)$. 
Furthermore by using Lemmas \ref{Demko} and \ref{cla2}, 
$(R^{(\ell)}_{n}G_{21n})\Lambda^{-1}_{1n}=O_P(K_n^{-2\ell+1}\vec{1}\vec{1}^\prime)$ can be obtained. 
By the repeat use of Lemma \ref{band} and Lemma \ref{Demko} in the same manner, we have 
\begin{eqnarray*}
R^{(\ell+1)}_{n}
&=&
O_P(K_n^{-2\ell+1}\vec{1}\vec{1}^\prime)G_{12n}\Lambda^{-1}_{2n}\\
&=&
O_P(K_n^{-2\ell-1}\vec{1}\vec{1}^\prime)\Lambda^{-1}_{2n}\\
&=&
O_P(K_n^{-2\ell}\vec{1}\vec{1}^\prime). \ \ \ \Box
\end{eqnarray*}

\vspace{5mm}

\begin{lemma}\label{eigen}
The maximum eigenvalues of $K_n(G_{jn}-G_j)$ and $K_nG_{12n}$ are asymptotically vanished.
\end{lemma}
\noindent
{\it Proof of Lemma \ref{eigen}}:\ Let $A_n=K_n(G_{1n}-G_1)=(a_{ij,n})_{ij}$. 
Then if $|i-j|\leq p$, $a_{ij,n}=o_P(1)$ and if $|i-j|>p$, $a_{ij,n}=0$ by Lemma \ref{G12}. 
Let $\lambda_{\max}(A_n)$ be the maximum eigenvalue of $A_n$. 
Then there exists $\vec{x}=(x_1\ \cdots\ x_{K_n+p})^\prime\in\mathbb{R}^{K_n+p}$ such that 
$$
A_n \vec{x}=\lambda_{\max}(A_n)\vec{x}.
$$
The $\vec{x}$ is eigenvector of $A_n$ belonging to $\lambda_{\max}(A_n)$.
Let $|x_m|$ be $\max\{|x_1|,\ \cdots,\ |x_{K_n+p}|\}$, we have 
\begin{eqnarray*}
|\lambda_{\max}(A_n)x_m|=\left|\sum_{j=1}^{K_n+p} a_{mj,n}x_j\right|. \label{eg}
\end{eqnarray*}
The $|\lambda_{\max}(A_n)|$ can be calculated as 
\begin{eqnarray*}
|\lambda_{\max}(A_n)|
&=&
\frac{1}{|x_m|}\left|\sum_{i=1}^{K_n+p} a_{mj,n}x_j\right|\\
&\leq &
\sum_{j=1}^{K_n+p} |a_{mj,n}|\frac{|x_j|}{|x_m|}\\
&\leq &
\sum_{j=1}^{K_n+p} |a_{mj,n}|\\
&=&
\sum_{j=m-p}^{m+p} |a_{mj,n}|\\
&=&
o_P(1)
\end{eqnarray*}
from the structure of $A_{n}$ .
The $K_nG_{12n}$ is also band matrix satisfying $K_nG_{12n}=o_P(\vec{1}\vec{1}^\prime)$. 
So we can prove that the maximum eigenvalue of $K_nG_{12n}$ is $o_P(1)$ by the same manner. $\Box$

\vspace{5mm}

\noindent
We are now in the position to give proofs of all results in Section 3.

\vspace{5mm}

\noindent
{\it Proof of Proposition \ref{initial}}:\ First we prove 
$|f_1^{(\ell)}(x_1)-f_{01}^{(\ell)}(x_1)|=O_P(K_n^{-2\ell})$. 
We have 
$$
f_1^{(\ell)}(x_1)-f_{01}^{(\ell)}(x_1)=-\vec{B}(x_1)^\prime (X_1^\prime X_1)^{-1}X_1^\prime(S_1S_2)^{\ell-1}S_1X_2\vec{b}_2^{(0)}
$$ 
and there exists $R^{(\ell-1)}_{n}=O_P(K_n^{-2(\ell-2)}\vec{1}\vec{1}^\prime)$ such that 
\begin{eqnarray*}
\vec{B}(x_1)^\prime (X_1^\prime X_1)^{-1}X_1^\prime(S_1S_2)^{\ell-1}S_1X_2\vec{b}_2^{(0)}
&=&\frac{1}{n}\vec{B}(x_1)^\prime R^{(\ell-1)}_{n}X_2^\prime S_1X_2\vec{b}_2^{(0)}\\
&=&\vec{B}(x_1)^\prime R^{(\ell-1)}_{n}G_{21n}\Lambda^{-1}_{jn}G_{12n}\vec{b}_2^{(0)}
\end{eqnarray*}
by Lemma \ref{S1S2}. 
We see from the proof of Lemma \ref{S1S2} that the $R^{(\ell)}_{n}$ consists of the product of $\Lambda_{1n}^{-1}$, $G_{12n}$, $\Lambda^{-1}_{2n}$ and $G_{21n}$ because 
\begin{eqnarray}
(S_1S_2)^{\ell}&=&(S_1S_2)(S_1S_2)\cdots (S_1S_2)\nonumber\\
&=&
X_1\{\Lambda_1^{-1}X_1^\prime X_2\Lambda_2^{-1}X_2^{\prime}(S_1S_2)\cdots (S_1S_2) X_1\Lambda_1^{-1}X_1^\prime X_2\Lambda_2^{-1}\}X_2^{\prime}\nonumber\\
&=&\frac{1}{n}X_1\{\Lambda_{1n}^{-1}G_{12n}\Lambda_{2n}^{-1}G_{21n}\cdots G_{21n}\Lambda_{1n}^{-1}G_{12n}\Lambda_{2n}^{-1}\}X_2^{\prime} \label{inte}.
\end{eqnarray}
Theorefore by Lemmas \ref{band}, \ref{Demko}, we have 
\begin{eqnarray*}
R^{(\ell-1)}_{n}G_{21n}\Lambda^{-1}_{jn}G_{12n}\vec{b}_2^{(0)}
&=&
R^{(\ell-1)}_{n}G_{21n}\Lambda^{-1}_{jn}O_P(K_n^{-2}\vec{1})\\
&=&
R^{(\ell-1)}_{n}G_{21n}O_P(K_n^{-1}\vec{1})\\
&=&
O_P(K_n^{-2\ell +1}\vec{1}),
\end{eqnarray*}
where $O_P(n^{\alpha}\vec{1})$ is the vector version of $O_P(n^\alpha \vec{1}\vec{1}^\prime)$. 
Because the $p+1$ components of $\vec{B}(x_1)$ are not 0 and others are 0 like the column of band matrix by property of $B$-spline basis, 
we have 
$$
\vec{B}(x_1)^\prime R^{(\ell-1)}_{n}G_{21n}\Lambda^{-1}_{jn}G_{12n}\vec{b}_2^{(0)}
=
O_P(K_n^{-2\ell+1})
$$ 
though the size of $\vec{B}(x_1)$ increases with $n$.   
Similarly, we see that  
$$
f_2^{(\ell)}(x_2)-f_{02}^{(\ell)}(x_2)=\vec{B}(x_2)^\prime \Lambda_2^{-1}X_2^\prime(S_1 S_2)^{\ell-1} S_1X_2\vec{b}_2^{(0)}
$$
becomes 
\begin{eqnarray*}
f_2^{(\ell)}(x_2)-f_{02}^{(\ell)}(x_2)
&=&
\vec{B}(x_2)^\prime \Lambda_{2n}^{-1}G_{21n}R^{(\ell-1)}_{n}G_{21n}\Lambda_{1n}^{-1}G_{12n}\vec{b}_2^{(0)}\\
&=&
O_P(K_n K_{n}^{-2}K_n^{-2(\ell -2)}K_n^{-2}K_n K_n^{-1})\\
&=&
o_P(K_n^{-2\ell+1}),
\end{eqnarray*}
which completes the proof. $\Box$

\vspace{5mm}

\noindent
{\it Proof of Proposition \ref{backfit}}:\ By Lemma \ref{S1S2}, we have 
\begin{eqnarray*}
f_{01}^{(\ell)}(x_1)
&=&
\vec{B}(x_1)^{\prime}(X_1^\prime X_1)^{-1}X_1^\prime\left\{I_n-\sum_{k=0}^{\ell-1} \{S_1S_2\}^k(I_n-S_1)\right\}\vec{y}\\
&=&
f_{01}^{(1)}(x_1)-\vec{B}(x_1)^{\prime}(X_1^\prime X_1)^{-1}X_1^\prime\sum_{k=1}^{\ell-1} \{S_1S_2\}^k(I_n-S_1)\vec{y}\\
&=&
f_{01}^{(1)}(x_1)-\vec{B}(x_1)^{\prime}\left\{ \sum_{k=1}^{\ell-1} R^{(k)}_n \right\}(X_2^\prime-G_{21n}\Lambda^{-1}_{1n}X_1^\prime)\frac{1}{n}\vec{y}
\end{eqnarray*}
and 
\begin{eqnarray*}
f_{02}^{(\ell)}(x_2)
&=& \vec{B}(x_2)^{\prime}\Lambda_2^{-1}X_2^\prime\sum_{k=0}^{\ell-1} \{S_1 S_2\}^k(I_n-S_1)\vec{y}\\
&=& f_{02}^{(1)}(x_2)+ \vec{B}(x_2)^{\prime}\Lambda_2^{-1}X_2^\prime\sum_{k=1}^{\ell-1} \{S_1 S_2\}^k(I_n-S_1)\vec{y}\\
&=& f_{02}^{(1)}(x_2)+ \vec{B}(x_2)^{\prime}\Lambda_{2n}^{-1}G_{21n}\left\{ \sum_{k=1}^{\ell-1} R^{(k)}_n \right\}(X_2^\prime-G_{21n}\Lambda^{-1}_{1n}X_1^\prime)\frac{1}{n}\vec{y}.
\end{eqnarray*}
We shall focus on the sum $\sum_{k=1}^{\ell-1} R^{(k)}_n$.
We put
$$
R_n=\lim_{\ell\rightarrow \infty}\sum_{k=1}^{\ell-1}R^{(k)}_n=\sum_{k=1}^{\infty}R^{(k)}_n=(r_{i,j,n})_{ij}.
$$
Then, since the backfitting algorithm converges for any $n$, $|r_{i,j,n}|$ is bounded for any $(i,j)$ and $n$.
And hence
$
r_{n} \equiv \max_{i,j}|r_{i,j,n}|
$
is also bounded for any $n$, which implies 
$$
R_{n}=O_P(\vec{1}\vec{1}).
$$

Let $\vec{f}=(f_1(x_{11})+f_2(x_{12})\ \cdots\ f_1(x_{n1})+f_2(x_{n2}))^\prime$. 
Then the absolute value of $h$-component of $n^{-1}X_j^\prime \vec{f}$ is
$$
\left|\frac{1}{n}\sum_{i=1}^n B_{-p+h}(x_{ij})(f_1(x_{i1})+f_2(x_{i2}))\right|\leq \max_{u,v\in(0,1)}\{|f_1(u)+f_2(v)|\}\frac{1}{n}\sum_{i=1}^n B_{-p+h}(x_{ij})=O_P(K_n^{-1})
$$
because $\max_{u,v\in(0,1)}\{|f_1(u)+f_2(v)|\}<\infty$. 
Hence $n^{-1}X_j^\prime \vec{f}=O_P(K_n^{-1}\vec{1})$ can be obtained. 
From (\ref{inte}) and the repeat use of Lemma \ref{band} and Lemma \ref{Demko}, we have 
$$
R_{n} n^{-1}X_j^\prime \vec{f}=O_P(K_n^{-1}\vec{1}).
$$
And direct calculation gives
\begin{eqnarray*}
|E[f_{01}^{(\infty)}(x_1)-f_{01}^{(1)}(x_1)]|
&=&
|\vec{B}(x_1)^{\prime}R_{n}(X_2^\prime-G_{21n}\Lambda^{-1}_{1n}X_1^\prime)\frac{1}{n}\vec{f}|\\
&=&
\vec{B}(x_1)^{\prime}O_P(K_n^{-1}\vec{1})+\vec{B}(x_1)^{\prime}O_P(K_n^{-2}\vec{1})\\
&=&O_P(K_n^{-1}\vec{1})
\end{eqnarray*}
because the $p+1$ components of $\vec{B}(x_1)$ are not 0 and others are 0. 
Here, $O_P(K_n^{-1}\vec{1})$ is the vector version of $O_P(K_n^{-1}\vec{1}\vec{1}^\prime)$. 
Similarly, 
we have 
\begin{eqnarray*}
|E[f_{02}^{(\infty)}(x_2)-f_{02}^{(1)}(x_2)]|
&=&
\left|\vec{B}(x_2)^{\prime}\Lambda_{2n}^{-1}G_{21n}R_{n}(X_2^\prime-G_{21n}\Lambda^{-1}_{1n}X_1^\prime)\frac{1}{n}\vec{f}\right|\\
&=&O_P(K_n K_n^{-2} 1 K_n^{-1})\\
&=&o_P(K_n^{-1}).
\end{eqnarray*}
Next, we consider 
$
V[f_{01}^{(\infty)}(x_1)-f_{01}^{(1)}(x_1)].
$
Let $\Sigma=\diag[\sigma^2(x_{11},x_{12})\ \cdots\ \sigma^2(x_{n1},x_{n2})]$. 
Then since $n^{-1}X_j^\prime \Sigma X_j=\Sigma_{j}+o_P(K_n^{-1}\vec{1}\vec{1}^\prime)=O_P(K_n^{-1}\vec{1}\vec{1}^\prime) (j=1,2)$,
we have 
\begin{eqnarray*}
&&V[f_{01}^{(\infty)}(x_1)-f_{01}^{(1)}(x_1)]\\
&&=R_{n}
\vec{B}(x_1)^{\prime}R_{n}(X_2^\prime-G_{21n}\Lambda^{-1}_{1n}X_1^\prime)\frac{1}{n}\Sigma\frac{1}{n}(X_2-X_1\Lambda^{-1}_{1n}G_{12n})R_{n}' \vec{B}(x_1)\\
&&=
\frac{1}{n}\vec{B}(x_1)^{\prime}R_{n}\Sigma_2R_{n}' \vec{B}(x_1)\\
&&\quad -2\frac{1}{n}\vec{B}(x_1)^{\prime}R_{n}\frac{1}{n}X_2^\prime \Sigma X_1\Lambda^{-1}_{1n}G_{12n}R_{n}' \vec{B}(x_1)\\
&&\quad \quad +\frac{1}{n}\vec{B}(x_1)^{\prime}R_{n}G_{21n}\Lambda^{-1}_{1n}\Sigma_1\Lambda^{-1}_{1n}G_{12n}R_{n}' \vec{B}(x_1)\\
&&=
\frac{1}{n}\vec{B}(x_1)^{\prime}R_{n}\Sigma_2R_{n}' \vec{B}(x_1)(1+o_P(1))\\
&&=
O_P((nK_n)^{-1}).
\end{eqnarray*}
Similarly, $V[f_{02}^{(\infty)}(x_2)-f_{02}^{(1)}(x_2)]$ can be calculated as 
\begin{eqnarray*}
&&V[f_{02}^{(\infty)}(x_2)-f_{02}^{(1)}(x_2)]\\
&&=
\frac{1}{n^2}\vec{B}(x_2)^{\prime}\Lambda_{2n}^{-1}G_{21n}R_{n}(X_2^\prime-G_{21n}\Lambda^{-1}_{1n}X_1^\prime)\\
&&\quad 
\times\Sigma (X_2-X_1\Lambda^{-1}_{1n}G_{12n})R_{n}'  G_{12n} \Lambda_{2n}^{-1}\vec{B}(x_2)\\
&&=
\frac{1}{n}\vec{B}(x_2)^{\prime}\Lambda_{2n}^{-1}G_{21n}R_{n} \Sigma_2R_{n}'  G_{12n} \Lambda_{2n}^{-1}\vec{B}(x_2)(1+o_P(1))\\
&&=
O_P(n^{-1}K_n K_n^{-2} 1 K_n^{-1} 1 K_n^{-2} K_n)\\
&&= o_P((nK_n)^{-1}).
\end{eqnarray*} 
Therefore, 
since 
$$
|f_{0j}^{(\infty)}(x_j)-f_{0j}^{(1)}(x_j)|=E[f_{0j}^{(\infty)}(x_j)-f_{0j}^{(1)}(x_j)]+O_P\left(\sqrt{V[f_{0j}^{(\infty)}(x_j)-f_{0j}^{(1)}(x_j)]}\right),
$$
we have 
$$
|f_{01}^{(\infty)}(x_1)-f_{01}^{(1)}(x_1)|=O_P(K_n^{-1})
$$
and 
$$
|f_{02}^{(\infty)}(x_2)-f_{02}^{(1)}(x_2)|=o_P(K_n^{-1}). \ \ \ \Box 
$$

\vspace{5mm}

\noindent
{\it Proof of Theorem \ref{mv}}:\ We see from Theorem 2 (a) of Claeskens et al. (2009) that  
\begin{eqnarray*}
E[f^{(1)}_{01}(x_1)]&=&f_1(x_1)+b_{1,\lambda}(x_1)+O_P(K_n^{-(p+1)})+o_P(\lambda_{1n}K_nn^{-1}),\\
V[f^{(1)}_{01}(x_1)]&=&V_1(x_1)+o_P(K_nn^{-1}).
\end{eqnarray*}
So we have 
\begin{eqnarray*}
E[f^{(\infty)}_{01}(x_1)]&=&E[f^{(1)}_{01}(x_1)]+O_P(K_n^{-1})\\
&=&f_1(x_1)+b_{1,\lambda}(x_1)+O_P(K_n^{-1})+o_P(\lambda_{1n}K_nn^{-1}),
\end{eqnarray*}
and 
\begin{eqnarray*}
V[f^{(\infty)}_{01}(x_1)]
&=&
V[f^{(1)}_{01}(x_1)]+2Cov\left(f_{01}^{(1)}(x_1),\theta_1\right)+V[\theta_1]\\
&=&
V_1(x_1)+o_P(K_nn^{-1}),
\end{eqnarray*}
where $\theta_{1}=f_{01}^{(\infty)}(x_1)-f_{01}^{(1)}(x_1)$.
This is in fact
\begin{eqnarray*}
\left|Cov\left(f_{01}^{(1)}(x_1),\theta_{1}\right)\right|
&\leq &\{V[f_{01}^{(1)}(x_1)]V[\theta_{1}]\}^{1/2}\\
&=&o_P(K_nn^{-1})
\end{eqnarray*}
and $V[\theta_1]=o_P(K_nn^{-1})$
from the proof of Proposition \ref{backfit}.
Furthermore, we also obtain 
\begin{eqnarray*}
E[f^{(\infty)}_{02}(x_2)]&=&f_2(x_2)+b_{2,\lambda}(x_2)+O_P(K_n^{-1})+o_P(\lambda_{2n}K_nn^{-1}),\\
V[f^{(\infty)}_{02}(x_2)]&=&V_2(x_2)+o_P(K_nn^{-1}).
\end{eqnarray*}
Finally, we calculate 
\begin{eqnarray*}
Cov(f^{(\infty)}_{01}(x_1),f^{(\infty)}_{02}(x_2))&=&Cov(f^{(1)}_{01}(x_1),f^{(1)}_{02}(x_2))+\sum_{j=1}^2 Cov\left(f_{0j}^{(1)}(x_j),\theta_{3-j}\right)+Cov (\theta_1,\theta_2), 
\end{eqnarray*}
where $\theta_{2}=f_{02}^{(\infty)}(x_2)-f_{02}^{(1)}(x_2)$.
Then we see that
\begin{eqnarray*}
Cov(f^{(1)}_{01}(x_1),f^{(1)}_{02}(x_2))
&=&
\frac{1}{n^2}\vec{B}(x_1)^\prime\Lambda_{1n}^{-1}X_1^\prime \Sigma X_2\Lambda_{2n}^{-1}\vec{B}(x_2)\\
&&+\frac{1}{n^2}\vec{B}(x_1)^\prime\Lambda_{1n}^{-1}X_1^\prime \Sigma X_1\Lambda_{1n}^{-1}G_{12n}\Lambda_{2n}^{-1}\vec{B}(x_2)\\
&=&
O(n^{-1} K_n K_n^{-2} K_n)\{ 1+o_{P}(1) \}\\
&=&O_P(n^{-1})
\end{eqnarray*}
because the absolute value of $(k,h)$-component of $n^{-1}X_1^\prime \Sigma X_2$ is 
\begin{eqnarray*}
\left|\frac{1}{n}\sum_{i=1}^n\sigma^2(x_{i1},x_{i2})B_{-p+k}(x_{i1})B_{-p+h}(x_{i2}) \right|
&\leq&
 \max_{u,v\in(0,1)}\{\sigma^2(u,v)\}\frac{1}{n}\sum_{i=1}^nB_{-p+k}(x_{i1})B_{-p+h}(x_{i2})\\
&=&
O_P(K_n^{-2}).
\end{eqnarray*}
In addition, for $j=1,2$, we find 
\begin{eqnarray*}
\left|Cov\left(f_{0j}^{(1)}(x_j),\theta_{3-j}\right)\right|
&\leq &
V[f_{0j}^{(1)}(x_j)]^{1/2}V[\theta_{3-j}]^{1/2}\\
&=&o_P(n^{-1})
\end{eqnarray*}
and $|Cov (\theta_1,\theta_2)|\leq \{ V[\theta_1]V[\theta_2] \}^{1/2}=o_P(n^{-1})$ from the proof of Proposition \ref{backfit}. 
This completes the proof. $\Box$

\vspace{5mm}

\noindent
{\it Proof of Theorem \ref{clt}}:\ If we prove 
\begin{eqnarray}
\left\{
V
\left[
\begin{array}{c}
\hat{f}_1(x_1)\\
\hat{f}_2(x_2)
\end{array}
\right]
\right\}^{-1/2}
\left[
\begin{array}{c}
\hat{f}_1(x_1)-E[\hat{f}_1(x_1)]\\
\hat{f}_2(x_2)-E[\hat{f}_2(x_2)]
\end{array}
\right]
\xrightarrow {D}N_2\left( 
\left[
\begin{array}{c}
0\\ 
0 
\end{array} 
\right]
,I_2\right), \label{norm}
\end{eqnarray}
then we have Theorem \ref{clt}.  
We rewrite $\hat{f}_j(x_j)$ as $\sum_{i=1}^n w_{i,jn}y_{i}$. 
For any $(a_1\ a_2)^\prime \in\mathbb{R}^2-\{\vec{0}\}$, we check 
$$
S_n\equiv a_1\hat{f}_1(x_1)+a_2\hat{f}_2(x_2)=\sum_{i=1}^n (a_1w_{i,1n}+a_2w_{i,2n})y_i
$$ 
satisfies the required Lyapunov condition.
First, we obtain 
\begin{eqnarray}
V[S_n]=a_1^2V[\hat{f}_1(x_1)]+2a_1a_2Cov(\hat{f}_1(x_1),\hat{f}_2(x_2))+a_2^2V[\hat{f}_2(x_2)]=O(K_nn^{-1}) \label{Sn}
\end{eqnarray}
by Theorem \ref{mv}. 
Next we note that the leading term of $w_{i,jn}$ 
$$
\vec{B}(x_j)^\prime \Lambda_{jn}^{-1}n^{-1}\vec{B}(x_{ij})\{ 1+o_{P}(1) \}\ \ j=1,2,
$$
because it is the $i$th component of $\vec{B}(x_j)^\prime \Lambda_{jn}^{-1}n^{-1}X_j^\prime$ with $X_j=(B_k(x_{ij}))_{ik}$.
Hence we have $w_{i,jn}=O_P(K_nn^{-1})$ and 
\begin{eqnarray}
&&E[|\{a_1w_{i,1n}+a_2w_{i,2n}\}y_i-\{a_1w_{i,1n}+a_2w_{i,2n}\}E[Y_i]|^{2+\delta}]\nonumber\\
&&= |a_1w_{i,1n}+a_2w_{i,2n}|^{2+\delta}E[|\varepsilon_i|^{2+\delta}]\nonumber\\
&&= O_P\left(\frac{K_n^{2+\delta}}{n^{2+\delta}}\right).\label{Ex}
\end{eqnarray} 
So it follows from (\ref{Sn}), (\ref{Ex}) and $K_n=O(n^\gamma)$ that 
\begin{eqnarray*}
\frac{1}{V[S_n]^{(2+\delta)/2}}\sum_{i=1}^n |a_1w_{i,1n}+a_2w_{i,2n}|^{2+\delta}E[|\varepsilon_i|^{2+\delta}]
&=&
O_P\left(K_n^{(2+\delta)/2}n^{-(2+\delta)/2}\right)O(n)O_P\left(\frac{K_n^{2+\delta}}{n^{2+\delta}}\right)\\
&=&
O_P\left(n^{\gamma(1+\delta/2)-\delta/2}\right).
\end{eqnarray*}
Therefore, for $\delta>2\gamma/(1-\gamma)$, 
$$
\frac{1}{V[S_n]^{(2+\delta)/2}}\sum_{i=1}^n E[|\{a_1w_{i,1n}+a_2w_{i,2n}\}y_i-\{a_1w_{i,1n}+a_2w_{i,2n}\}E[Y_i]|^{2+\delta}]
=o_P(1).
$$
By Lyapunov theorem and Cram\'{e}r-Wold Device, we get (\ref{norm}). 

Consequently, since asymptotic bias of $\hat{f}_j(x_j)$ is $b_{j,\lambda}(x_j)=O_P(\lambda_{jn}K_nn^{-1})=o_P(K_n^{1/2}n^{-1/2})$, 
Theorem \ref{clt} has been obtained.  $\Box$

\vspace{5mm}

\noindent
{\it Proof of Theorem 3}:\ We show that the $H(L)$ becomes positive definite as $n\rightarrow \infty$. 
Now, $H(L)$ is divided into
\begin{eqnarray*}
H(L)&=&
\left[
\begin{array}{cc}
X_1^\prime X_1 & X_1^\prime X_2\\
X_2^\prime X_1&X_2^\prime X_2
\end{array}
\right]
+
\left[
\begin{array}{cc}
\lambda_{1n}Q_m & O\\
O&\lambda_{2n}Q_m
\end{array}
\right]\\
&\equiv& H_1 + H_2.
\end{eqnarray*}
Then it is known that $Q_m$ has eigenvalue 0, 
hence the $H_2$ is nonnegative definite. 
We show that the $H_1$ becomes positive definite as $n\rightarrow \infty$. 
For any $\vec{z}\in\mathbb{R}^{2(K_n+p)}$ with $\vec{z}^\prime\vec{z}=1$, 
\begin{eqnarray*}
\frac{K_n}{n}\vec{z}^\prime H_1\vec{z}
&=&
\vec{z}^\prime\left[
\begin{array}{cc}
K_nG_1 & O\\
O&K_nG_2
\end{array}
\right]\vec{z}
+
\vec{z}^\prime\left[
\begin{array}{cc}
K_n(G_{1n}-G_1) & K_nG_{12n}\\
K_nG_{21n}&K_n(G_{2n}-G_2)
\end{array}
\right]\vec{z}\\
&\equiv&  \vec{z}^\prime H_{11}\vec{z} + \vec{z}^\prime H_{12}\vec{z}.
\end{eqnarray*} 
By (6.10) of Agarwal and Studden (1980), 
we can find $\vec{z}^\prime H_{11}\vec{z}>0$. 
Now we show $\vec{z}^\prime H_{12}\vec{z}=o_P(1)$. 
We write $\vec{z}=(\vec{z}_1^\prime\ \vec{z}_2^\prime)^\prime$, where $\vec{z}_j\in\mathbb{R}^{K_n+p}$. 
Then since $\vec{z}^\prime \vec{z}=\vec{z}_1^\prime \vec{z}_1+\vec{z}_2^\prime \vec{z}_2=1$, $\vec{z}_j^\prime \vec{z}_j\leq 1 (j=1,2)$.
So we get 
\begin{eqnarray*}
\vec{z}^\prime H_{12}\vec{z}
&=&
\vec{z}_1^\prime K_n(G_{1n}-G_1)\vec{z}_1 +2\vec{z}_1^\prime K_nG_{12n}\vec{z}_2+\vec{z}_2^\prime K_n(G_{2n}-G_2)\vec{z}_2.
\end{eqnarray*}  
The maximum eigenvalue of $K_n(G_{jn}-G_j)$ is $o_p(1)$ from Lemma \ref{eigen}, 
we have 
\begin{eqnarray*}
\vec{z}_j^\prime K_n(G_{jn}-G_j)\vec{z}_j=\vec{z}_j^\prime\vec{z}_j \frac{\vec{z}_j^\prime K_n(G_{jn}-G_j)\vec{z}_j}{\vec{z}_j^\prime\vec{z}_j}=o_P(1),\ \ j=1,2.
\end{eqnarray*}
Similarly, we get
\begin{eqnarray*}
|\vec{z}_1^\prime K_nG_{12n}\vec{z}_2|&\leq& \sqrt{\vec{z}_1^\prime\vec{z}_1}\sqrt{\vec{z}_2^\prime\vec{z}_2\frac{\vec{z}_2^\prime K_nG_{12n}\vec{z}_2}{\vec{z}_2^\prime\vec{z}_2}}\\
&=&\sqrt{\vec{z}_1^\prime\vec{z}_1}\sqrt{\vec{z}_2^\prime\vec{z}_2o_P(1)}\\
&=&o_P(1)
\end{eqnarray*}
because the maximum eigenvalue of $K_nG_{12n}$ is also $o_p(1)$ as shown in Lemma \ref{eigen}. 
Above evaluations are combined into $\vec{z}^\prime H_{12}\vec{z}=o_P(1)$. 
Consequently, we obtain
\begin{eqnarray*}
\vec{z}^\prime H_1\vec{z}
&=&
\frac{n}{K_n}\left(\vec{z}^\prime H_{11}\vec{z} + \vec{z}^\prime H_{12}\vec{z}\right)\\
&=&
\frac{n}{K_n} \vec{z}^\prime H_{11}\vec{z}\left(1+\frac{\vec{z}^\prime H_{12}\vec{z}}{\vec{z}^\prime H_{11}\vec{z}}\right)\\
&=&
\frac{n}{K_n} \vec{z}^\prime H_{11}\vec{z}(1+o_P(1))\\
&>&0. \ \ \ \Box
\end{eqnarray*}

\def\bibname{Reference}

\end{document}